\documentclass[11pt]{article}
\usepackage{amsmath,amsfonts}
\usepackage{verbatim}
\usepackage{latexsym}
\usepackage{graphicx}
\usepackage{float}
\usepackage{color}
\usepackage{cite}
\usepackage{subfigure}
\allowdisplaybreaks
\makeatletter
\@addtoreset{equation}{section}

\begin{document}

\newcommand{\E}{\mathbb{E}}
\newcommand{\PP}{\mathbb{P}}
\newcommand{\RR}{\mathbb{R}}

\newtheorem{theorem}{Theorem}[section]
\newtheorem{lemma}[theorem]{Lemma}
\newtheorem{coro}[theorem]{Corollary}
\newtheorem{defn}[theorem]{Definition}
\newtheorem{assp}[theorem]{Assumption}
\newtheorem{expl}[theorem]{Example}
\newtheorem{prop}[theorem]{Proposition}
\newtheorem{rmk}[theorem]{Remark}

\newcommand\tq{{\scriptstyle{3\over 4 }\scriptstyle}}
\newcommand\qua{{\scriptstyle{1\over 4 }\scriptstyle}}
\newcommand\hf{{\textstyle{1\over 2 }\displaystyle}}
\newcommand\hhf{{\scriptstyle{1\over 2 }\scriptstyle}}

\newcommand{\proof}{\noindent {\it Proof}. }
\newcommand{\eproof}{\hfill $\Box$} 

\def\a{\alpha} \def\g{\gamma}
\def\e{\varepsilon} \def\z{\zeta} \def\y{\eta} \def\o{\theta}
\def\vo{\vartheta} \def\k{\kappa} \def\l{\lambda} \def\m{\mu} \def\n{\nu}
\def\x{\xi}  \def\r{\rho} \def\s{\sigma}
\def\p{\phi} \def\f{\varphi}   \def\w{\omega}
\def\q{\surd} \def\i{\bot} \def\h{\forall} \def\j{\emptyset}

\def\be{\beta} \def\de{\delta} \def\up{\upsilon} \def\eq{\equiv}
\def\ve{\vee} \def\we{\wedge}

\def\F{{\cal F}}
\def\T{\tau} \def\G{\Gamma}  \def\D{\Delta} \def\O{\Theta} \def\L{\Lambda}
\def\X{\Xi} \def\S{\Sigma} \def\W{\Omega}
\def\M{\partial} \def\N{\nabla} \def\Ex{\exists} \def\K{\times}
\def\V{\bigvee} \def\U{\bigwedge}

\def\1{\oslash} \def\2{\oplus} \def\3{\otimes} \def\4{\ominus}
\def\5{\circ} \def\6{\odot} \def\7{\backslash} \def\8{\infty}
\def\9{\bigcap} \def\0{\bigcup} \def\+{\pm} \def\-{\mp}
\def\la{\langle} \def\ra{\rangle}

\def\tl{\tilde}
\def\trace{\hbox{\rm trace}}
\def\diag{\hbox{\rm diag}}
\def\for{\quad\hbox{for }}
\def\refer{\hangindent=0.3in\hangafter=1}

\newcommand\wD{\widehat{\D}}
\title{
\bf Numerical scheme for delay-type stochastic McKean-Vlasov equations driven by fractional Brownian motion}

\author{
{\bf Shuaibin Gao${}^{1}$, Qian Guo$^{1}$, Zhuoqi Liu$^{1}$, Chenggui Yuan$^{2}$\thanks{The corresponding author. Email: C.Yuan@swansea.ac.uk } }
\\
${}^1$ Department of Mathematics, 
Shanghai Normal University, \\
Shanghai, 200234, China. \\
${}^2$ Department of Mathematics, 
Swansea University, Bay Campus,\\
 Swansea, SA1 8EN, UK. \\
 }

\date{}

\maketitle

\begin{abstract}
This paper focuses on the numerical scheme for 
delay-type stochastic McKean-Vlasov equations (DSMVEs) driven by fractional Brownian motion with Hurst parameter $H\in (0,1/2)\cup (1/2,1)$.
The existence and uniqueness of the solutions to such DSMVEs whose drift coefficients contain polynomial delay terms are proved by exploting the Banach fixed point theorem. Then the propagation of chaos between interacting particle system and non-interacting system in $\mathcal{L}^p$ sense is shown. 
We find that even if the delay term satisfies the polynomial growth condition, the unmodified classical  Euler-Maruyama scheme still can approximate the corresponding interacting particle system without the particle corruption. The convergence rates are revealed for $H\in (0,1/2)\cup (1/2,1)$.
Finally, as an example that closely fits the original equation, a stochastic opinion dynamics model with both extrinsic memory and intrinsic memory is simulated to illustrate the plausibility of the theoretical result.

 \par \noindent

\end{abstract}

\section{Introduction}\label{cccsec01}

As is known to all, the fractional Brownian motion (fBm) with Hurst parameter $H\in (0,1)$ is a natural generalization of the usual Brownian motion.
In light of the extrinsic memory impact of fBm on system, the theories about stochastic differential equations (SDEs) driven by fBm 
have been systematacially studied in \cite{ccc1,ccc2}.  
Since the analytical solution to SDE  driven by fBm cannot be expressed explicitly in many scenarios, the investigation of numerical schemes becomes crucial. Many types of SDEs  driven by fBm can be approximated by backward Euler-Maruyama (EM) scheme \cite{ccc3,ccc4,ccc5}, $\theta$-EM scheme \cite{ccc6,ccc7}, Milstein-type scheme \cite{ccc8,ccc9}, Crank-Nicolson scheme \cite{ccc10}, tamed EM scheme \cite{ccc18}, truncated EM sheme \cite{ccc11} and so on.

When the coefficients of SDEs are related to the laws of state variables, the equations are called stochastic McKean-Vlasov equations (SMVEs), also called mean-field SDEs or distribution-dependent SDEs.
For SMVEs driven by fBm with Hurst parameter $H\in (0,1/2)\cup (1/2,1)$, the wellposedness and the Bismut formula for the Lions derivative were presented in \cite{ccc12}. Then EM scheme was exploited to approximate such SMVEs driven by fBm in \cite{ccc34}.
In reality, there are already many papers analyzing the numerical schemes of SMVEs driven by standard Brownian motion (fBm with  $H=1/2$) via stochastic particle method, such as implicit EM scheme and tamed EM scheme in \cite{ccc35},  tamed Milstein scheme in \cite{ccc37,ccc38},
multi-level Monte-Carlo
scheme in \cite{ccc36},
adaptive EM and Milstein scheme in \cite{ccc49}, and split‐step EM scheme in \cite{ccc48}.

When the influence of the intrinsic memory on system is taken into account, the delay-type SDEs driven by  fBm have been investigated in \cite{ccc39,ccc40,ccc41,ccc42,ccc43,ccc44,ccc45}.
However, there are few work on the delay-type stochastic McKean-Vlasov equations (DSMVEs) driven by  fBm yet. 
To fill this gap,
this paper focuses on a class of DSMVEs driven by  fBm with $H\in (0,1/2)\cup (1/2,1)$ and their numerical schemes via  stochastic particle method. 
We first give the existence and uniqueness of solution to DSMVE  driven by  fBm by exploiting the Banach fixed point theorem.
Then the propagation of chaos in $\mathcal{L}^p$ sense between non-interacting particle system and interacting particle system is presented so that the unmodified classical EM schemes can be established.
The convergence rates of classical EM schemes are given   for $H\in (0,1/2)\cup (1/2,1)$.

Undoubtedly, the unmodified classical EM scheme fails to approximate the  SDEs driven by standard Brownian motion if the coefficients are superlinear \cite{ccc13}.
However, what needs to be emphasized in our paper is that the present state variable in the drift coefficient satisfies the global Lipschitz condition but the past state variable can grow polynomially, and the unmodified classical EM scheme still works in this case, which means that the particle corruption shown in \cite{ccc35} does not occur.

 Stochastic opinion dynamics model (SODM), 
which can reflect changes in people's online or offline opinions in the social scenarios, is an crucial tool to formulate promotional plans and information campaigns, so SODM has been extensively studied \cite{ccc14,ccc15,ccc17}.
Compared to the standard Brownian motion, SODM driven by  fBm better reflects the influence of extrinsic memory due to the long-range (or short-range ) dependence of fBm.
Additionally, the impact of delay state variable and its distribution information on SODM should certainly be taken into account, as communications could be delayed in both online and offline between individuals in real-life scenarios \cite{ccc16}.
We will find that our target DSMVEs driven by fBm can well match the SODM with both extrinsic memory and intrinsic memory. Therefore, the numerical analysis of DSMVEs driven by fBm is an important task for further research on SODM.

The rest of the paper is structured as follows. In Section \ref{cccsec02}, some notations and  important lemmas are introduced, then we give the wellposedness of DSMVE driven by  fBm. Section \ref{cccsec03} aims to reveal the propagation of chaos in $\mathcal{L}^p$ sense. 
For $H\in (0,1/2)\cup (1/2,1)$, the convergence rates of classical EM scheme for interacting particle system are shown in Section \ref{cccsec04}.
As a suitable example, a SODM is simulated in Section \ref{cccsec05}.

\section{Preliminaries}\label{cccsec02}

Let $|\cdot|$ be the Euclidean norm for vector and the trace norm for matrix.
Denote $a_1\vee a_2=\max\{a_1,a_2\}$, $a_1\wedge a_2=\min\{a_1,a_2\}$ for real numbers $a_1,a_2$.
Let $B^H=\{B_t^H, t\geq 0\}$ be a fBm with Hurst parameter $H\in (0,1)$ defined
on the probability space $(\Omega, \mathcal{F}, \mathbb{P})$, i.e., $B^H$ is
a centered Gaussian process  with covariance function 
$R_H(t,s)=\frac{1}{2}(t^{2H}+s^{2H}-|t-s|^{2H})$ for any $s,t\geq 0$. 
Let $\mathbb{I}_S$ be the indicator function for a set $S$, i.e., $\mathbb{I}_S(x)=1$ if $x \in S$, otherwise, $\mathbb{I}_S(x)=0$.
Let $\lfloor a\rfloor$ be the largest integer which does not exceed $a$.
Denote by $\mathcal{L}^p =\mathcal{L}^p( \Omega, \mathcal{F},  \mathbb{P} \big )$  the set of random variables $X$ with  expectation $\mathbb{E}|X|^p<\infty$ for $p\geq1$.
Let $\mathcal{C}=\mathcal{C}([-\rho,0];\mathbb{R}^d)$ be the family of all continuous functions $\varphi$ from $[-\rho,0] $ to $\mathbb{R}^d$ with the norm $\|\varphi\|=\text{sup}_{-\rho \le \theta \le 0}|\varphi (\theta)|$.
For a positive integer $N$, let $\mathbb{S}_N=\{1,2,\cdots, N\}$.

Denote by  $\delta_{y}(\cdot)$  the Dirac measure at point $y\in \mathbb{R}^{d}$.
Denote by $\mathcal{P}(\mathbb{R}^{d})$ the family of all probability measures on $\mathbb{R}^{d}$. For $q\geq1$, define
$$\mathcal{P}_{q}(\mathbb{R}^{d})=\left\{\mu\in\mathcal{P}(\mathbb{R}^{d}):\left(\int_{\mathbb{R}^{d}}|y|^{q}\mu(dy)\right)^{1/q}<\infty\right\},$$
and set $\mathcal{W}_{q}(\mu)=\left(\int_{\mathbb{R}^{d}}|y|^{q}\mu(dy)\right)^{1/q}$ for any $\mu\in \mathcal{P}_{q}(\mathbb{R}^{d})$.
For $q\geq1$, let $\mathcal{C}([-\rho,T];\mathcal{P}_{q}(\mathbb{R}^{d}))$ be the family of all continuous measures $\mu$ from $[-\rho,T] $ to $\mathcal{P}_{q}(\mathbb{R}^{d})$.
For $q\geq1$,  the Wasserstein distance of $\mu,\nu\in \mathcal{P}_{q}(\mathbb{R}^{d})$ is defined by
$$\mathbb{W}_{q}(\mu,\nu)=\inf_{\pi\in \mathfrak{C}(\mu,\nu)}\left(\int_{\mathbb{R}^{d}\times\mathbb{R}^{d}}|y_1-y_2|^{q}\pi(dy_1,dy_2)\right)^{1/q},$$
where $\mathfrak{C}(\mu,\nu)$ is the family of all couplings for $\mu,\nu$.

In this paper, consider the following DSMVE driven by  fBm 
\begin{equation}\label{ccc1.1}
d\Upsilon(t)
=\alpha_{t}\left(\Upsilon(t),\Upsilon(t-\rho),\mathbb{L}_{\Upsilon(t)},\mathbb{L}_{\Upsilon(t-\rho)}\right)dt+\beta_{t}\left(\mathbb{L}_{\Upsilon(t)},\mathbb{L}_{\Upsilon(t-\rho)}\right) dB_t^H,~~~ t\in[0,T],
\end{equation}
with the initial value $\{\Upsilon (\theta): -\rho\leq \theta \leq 0\}=\xi$, which is an $\mathcal{F}_0$-measurable $\mathcal{C}$-valued
random variable with $\mathbb{E}\|\xi\|^{\check{p}}<\infty$ for any $\check{p}>0$.
Here, $\mathbb{L}_{\Upsilon(\cdot)}$ is the distribution of  $\Upsilon(\cdot)$. Moreover, $\alpha:[0,T]\times\mathbb{R}^{d}\times\mathbb{R}^{d}\times\mathcal{P}_{2}(\mathbb{R}^{d})\times\mathcal{P}_{2}(\mathbb{R}^{d})\rightarrow\mathbb{R}^{d},$ $\beta:[0,T]\times\mathcal{P}_{2}(\mathbb{R}^{d})\times\mathcal{P}_{2}(\mathbb{R}^{d})\rightarrow\mathbb{R}^{d\times d}$ are Borel measurable functions, and $B_t^H$ is a $d$-dimensional fBm with $H\in (0,1/2)\cup (1/2,1)$. 
Note that $\int_0^T \beta_{t}\left(\mathbb{L}_{\Upsilon(t)},\mathbb{L}_{\Upsilon(t-\rho)}\right) dB_t^H$
is treated as a Wiener integral with respect to fBm since the diffusion coefficient is a deterministic function.

To get the wellposedness of solution to (\ref{ccc1.1}) by exploiting the Banach fixed point theorem,
we consider the distribution-independent delay-type SDE driven by fBm of the form
\begin{equation} \label{cccada1}
	d\hat{\Upsilon}(t)
	=\hat{\alpha}_{t}\left(\hat{\Upsilon}(t),\hat{\Upsilon}(t-\rho)\right)dt+\hat{\beta}_{t} dB_t^H,~~~ t\in[0,T],
\end{equation}
with the initial value $\xi$.
Moreover, $\hat{\alpha}:[0,T]\times\mathbb{R}^{d}\times\mathbb{R}^{d}\rightarrow\mathbb{R}^{d},$ and $\hat{\beta}:[0,T]\rightarrow\mathbb{R}^{d}$ are Borel measurable functions. Assume that $\hat{\alpha}_{t}(0,0)$ and $\hat{\beta}_{t}$ are bounded for any $t\in[0,T]$.
We give the following important maximal inequality about fBm for the case $H\in (0,1/2)\cup (1/2,1)$, which has been proved by Theorem 1.2 in \cite{ccc30} and Theorem 2.1 in \cite{ccc31}.
\begin{lemma}\label{cccineq1}
	For any $\check{p}>0$ and $H\in (0,1/2)\cup (1/2,1)$,  there exist two  constants $c_{H,\check{p}},C_{H,\check{p}}>0$ such that
	\begin{equation*}
		c_{H,\check{p}}\mathbb{E}(\tau^{\check{p}H})\leq\mathbb{E}\big(\sup_{t\in[0,\tau]}|B_t^H|^{\check{p}}\big)\leq C_{H,\check{p}}\mathbb{E}(\tau^{\check{p}H}),	
	\end{equation*}
for any stopping time $\tau$ of $B_t^H$.
\end{lemma}

\begin{assp}\label{cccass3}
There exist two constants $\bar{K}_1>0$ and $l\geq 1$ such that
	\begin{equation*}
		|\hat{\alpha}_{t}(x_1,x_2)-\hat{\alpha}_{t}(y_1,y_2)|\leq \bar{K}_1[|x_1-y_1|+(1+|x_2|^l+|y_2|^l)|x_2-y_2|],
	\end{equation*}
 for any $t\in [0,T]$ and $x_1,x_2,y_1,y_2\in\mathbb{R}^d$.
\end{assp}
By Assumption \ref{cccass3}, one can see that there exists a $\bar{C}_{1}>0$ such that
\begin{equation*}
	\left|\hat{\alpha}_{t}(x_1,x_2)\right|\leq \bar{C}_{1}\left(1+|x_1|+|x_{2}|^{l+1}\right),
\end{equation*}
for any $t\in [0,T]$ and $x_1,x_2\in\mathbb{R}^d$.
The following lemma reveals that there exists a uinque global solution to (\ref{cccada1}) under Assumption \ref{cccass3}.

\begin{lemma}\label{cccsdeexistence}
Let Assumption \ref{cccass3} hold and $H\in (0,1/2)\cup (1/2,1)$. Then there exists a unique global solution $\hat{\Upsilon}(t)$ to (\ref{cccada1}), and it satisfies, for any $\bar{p}\geq2$,
\begin{equation*}
\mathbb{E}\big(\sup_{t\in[0,T]}|\hat{\Upsilon}(t)|^{\bar{p}}\big)\leq C_{\bar{p},T,H,\bar{K}_1,\|\xi\|,l,\bar{p}_m}.
\end{equation*}
\end{lemma}

\noindent
{\it Proof}.
	In view of \cite{ccc20}, one can see that (\ref{cccada1}) admits a unique local solution when the drift coefficient is local Lipschitz continuous and the diffusion coefficient is a function of $t$. So, to achieve the goal of this lemma, we just need to prove 
	
	\begin{equation*}
		\mathbb{E}\big(\sup_{t\in[0,T]}|\hat{\Upsilon}(t)|^{\bar{p}}\big)\leq C,~~~\text{for any}~ \bar{p}\geq2.
	\end{equation*}
From (\ref{cccada1}), using H\"older's inequality, Assumption \ref{cccass3} and Lemma \ref{cccineq1} leads to 
	\begin{equation*}
\begin{split}
&\mathbb{E}\big(\sup_{t\in[0,T]}|\hat{\Upsilon}(t)|^{\bar{p}}\big)\\
	\leq& C_{\bar{p},T,H,\|\xi\|}+C_{p,T}\mathbb{E}\int_{0}^{T}\big|\alpha_{s}(\hat{\Upsilon}(s),\hat{\Upsilon}(s-\rho))\big|^{\bar{p}}ds+\mathbb{E}\left(\sup_{t\in[0,T]}|\int_0^t\hat{\beta}_{s} dB_s^H|^{\bar{p}}\right)\\
		\leq& C_{\bar{p},T,H,\|\xi\|}+C_{p,T,\bar{K}_1}\mathbb{E}\int_{0}^{T}\big[1+|\hat{\Upsilon}(s)|^{\bar{p}}+|\hat{\Upsilon}(s-\rho)|^{\bar{p}(l+1)}\big]ds\\
			\leq& C_{\bar{p},T,H,\bar{K}_1,\|\xi\|}\left[1+\int_{0}^{T}\mathbb{E}\big(\sup_{s\in[0,t]}|\hat{\Upsilon}(s)|^{\bar{p}}\big)dt+\int_{0}^{T}\mathbb{E}|\hat{\Upsilon}(s-\rho)|^{\bar{p}(l+1)}ds\right].
\end{split}
	\end{equation*}
The Gronwall inequality means that
		\begin{equation}\label{ccc33}
		\begin{split}
			\mathbb{E}\big(\sup_{t\in[0,T]}|\hat{\Upsilon}(t)|^{\bar{p}}\big)
			\leq C_{\bar{p},T,H,\bar{K}_1,\|\xi\|}\big[1+\mathbb{E}\big(\sup_{t\in[0,T]}|\hat{\Upsilon}(t-\rho)|^{\bar{p}(l+1)}\big)\big].
		\end{split}
	\end{equation}
Define a sequence as
\begin{equation*}
	\bar{p}_{m}=(2-m+\lfloor\frac{T}{\rho}\rfloor)\bar{p}(l+1)^{1-m+\lfloor\frac{T}{\rho}\rfloor},~~~m=1,2 \cdots,\lfloor\frac{T}{\rho}\rfloor+1.
\end{equation*}
Obviously, it holds that $\bar{p}_{\lfloor\frac{T}{\rho}\rfloor+1}=\bar{p}$ and $\bar{p}_{m+1}(l+1)<\bar{p}_{m}$ for $m=1,2 \cdots,\lfloor\frac{T}{\rho}\rfloor$.
For $t\in[0,\rho]$, we get from (\ref{ccc33}) that
	\begin{equation*}
	\begin{split}
		\mathbb{E}\big(\sup_{t\in[0,\rho]}|\hat{\Upsilon}(t)|^{\bar{p}_1}\big)
		\leq& C_{\bar{p},T,H,\bar{K}_1,\|\xi\|}\big[1+\mathbb{E}\big(\sup_{t\in[0,\rho]}|\hat{\Upsilon}(t-\rho)|^{{\bar{p}_{1}}(l+1)}\big)\big]\\
		\leq&C_{\bar{p},T,H,\bar{K}_1,\|\xi\|,l,\bar{p}_1}.
	\end{split}
\end{equation*}
Next, for $t\in[0,2\rho]$, using H\"older's inequality, (\ref{ccc33}) and $\bar{p}_{m+1}(l+1)<\bar{p}_{m}$ leads to 
	\begin{equation*}
	\begin{split}
		\mathbb{E}\big(\sup_{t\in[0,2\rho]}|\hat{\Upsilon}(t)|^{\bar{p}_2}\big)
		\leq& C_{\bar{p},T,H,\bar{K}_1,\|\xi\|}\Big[1+\left[\mathbb{E}\big(\sup_{t\in[0,2\rho]}|\hat{\Upsilon}(t-\rho)|^{{\bar{p}_{1}}}\big)\right]^{\frac{\bar{p}_2(l+1)}{\bar{p}_1}}\Big]\\
		\leq&C_{\bar{p},T,H,\bar{K}_1,\|\xi\|,l,\bar{p}_1,\bar{p}_2}.
	\end{split}
\end{equation*}
Repeating this procedure gives 
	\begin{equation*}
		\mathbb{E}\left(\sup_{t\in\left[0,(\lfloor\frac{T}{\rho}\rfloor+1)\rho\right]}|\hat{\Upsilon}(t)|^{\bar{p}}\right)\leq C_{\bar{p},T,H,\bar{K}_1,\|\xi\|,l,\bar{p}_m},
	\end{equation*}
	where $m=1,2 \cdots,\lfloor\frac{T}{\rho}\rfloor+1$.
	\eproof
	
When $H=1/2$, $B_t^H$ becomes a standard Brownian motion, and the results have been presented in \cite{ccc24}.
Then under the following condition, the wellposedness of (\ref{ccc1.1}) is given by means of Lemma \ref{cccsdeexistence}.
	\begin{assp}\label{cccass1}
		There exist three constants $K_1,K_2>0$ and $l\geq 1$ such that
		\begin{equation*}
			|\alpha_{t}(x_1,x_2,\mu_1,\mu_2)-\alpha_{t}(y_1,y_2,\nu_1,\nu_2)|\leq K_1[|x_1-y_1|+(1+|x_2|^l+|y_2|^l)|x_2-y_2|+\mathbb{W}_{2}(\mu _{1},\nu_{1} )+\mathbb{W}_{2}(\mu _{2},\nu_{2} )],
		\end{equation*}
		\begin{equation*}
			|\beta_{t}(\mu_1,\mu_2)-\beta_{t}(\nu_1,\nu_2)|\leq K_2(\mathbb{W}_{2}(\mu _{1},\nu_{1} )+\mathbb{W}_{2}(\mu _{2},\nu_{2} )),
		\end{equation*}
	 for any $t\in [0,T]$,  $x_1,x_2,y_1,y_2\in\mathbb{R}^d$ and $\mu_1,\mu_2,\nu_1,\nu_2\in \mathcal{P}_2(\mathbb{R}^d)$.
	\end{assp}

\begin{lemma}\label{ccclmm1}
	For any $\mu\in \mathcal{P}_{2}(\mathbb{R}^{d})$, we have $\mathbb{W}_{2}(\mu,\delta_{0})=\mathcal{W}_{2}(\mu)$.
\end{lemma}
The proof of above lemma can be found in
Lemma 2.3 of \cite{ccc21}.
By Assumption \ref{cccass1} and Lemma \ref{ccclmm1}, one can see that there exist constants $\bar{C}_{2},\bar{C}_{3}>0$ such that
\begin{equation*}
	\left|{\alpha}_{t}(x_1,x_2,\mu_1,\mu_2)\right|\leq \bar{C}_{2}\left(1+|x_1|+|x_{2}|^{l+1}+\mathcal{W}_{2}(\mu _{1} )+\mathcal{W}_{2}(\mu _{2} )\right),
\end{equation*}
	\begin{equation*}
	|\beta_{t}(\mu_1,\mu_2)|\leq \bar{C}_{3}(1+\mathcal{W}_{2}(\mu _{1} )+\mathcal{W}_{2}(\mu _{2} )),
\end{equation*}
for any $t\in [0,T]$, $x_1,x_2\in\mathbb{R}^d$ and $\mu_1,\mu_2\in \mathcal{P}_2(\mathbb{R}^d)$.
By borrowing the proof ideas of Theorem 3.3 in \cite{ccc21}, Theorem 1.4 in \cite{ccc47} and Theorem 1.12.1 in \cite{ccc1}, we give the following theorem. 
\begin{theorem}\label{cccexactbound}
If $H\in (1/2,1)$, let Assumption \ref{cccass1} hold. If $H\in (0,1/2)$, let $\alpha$ satisfy Assumption \ref{cccass1} and $\beta$ only depend on $t$ (i.e., $\beta$ does not depend on the distribution).
Then DSMVE driven by fBm (\ref{ccc1.1}) admits a unique global solution $\Upsilon(t)$ satisfying for any $\bar{p}\geq2$ and $T>0$,
\begin{equation*}
	\mathbb{E}\big(\sup_{t\in[0,T]}|\Upsilon(t)|^{\bar{p}}\big)\leq C_{\bar{p},T,H,K_1,K_2,\|\xi\|,\bar{p}_m}.
\end{equation*}
\end{theorem}
\noindent
{\it Proof}.
We first show the assertion for $H\in (1/2,1)$.
For $x,y\in\mathbb{R}^d$ and $\mu_{.} \in\mathcal{C}\left([-\rho,T];\mathcal{P}_2(\mathbb{R}^d)\right)$, let $\alpha_{t,\rho}^{\mu}(x,y)=\alpha_{t}(x,y,\mu_{t},\mu_{t-\rho})$ and $\beta_{t,\rho}^{\mu}=\beta_{t}(\mu_{t},\mu_{t-\rho})$.
Consider the auxiliary SDE of the form
\begin{equation}\label{cxcxa1}
	d\Upsilon^{\mu}(t)
	=\alpha_{t,\rho}^{\mu}(\Upsilon^{\mu}(t),\Upsilon^{\mu}\left(t-\rho)\right)dt+\beta_{t,\rho}^{\mu} dB_t^H,~~~t\in[0,T],
\end{equation}
with the initial value $\Upsilon_0^\mu=\xi$.
Under Assumption \ref{cccass1}, we derive from Lemma \ref{cccsdeexistence} that  SDE (\ref{cxcxa1}) admits a unique global solution in a strong sense and it satisfies
\begin{equation}\label{cxcxa2}
	\mathbb{E}\big(\sup_{t\in[0,T]}|\Upsilon^{\mu}(t)|^{\bar{p}}\big)\leq C_{\bar{p},T,H,\bar{K}_1,\|\xi\|,l,\bar{p}_m},
\end{equation}
for any $\bar{p}\geq2$.
Define an operator
\begin{equation*}
	\Phi_t:\mathcal{C}\left([-\rho,T];\mathcal{P}_2(\mathbb{R}^d)\right)\rightarrow\mathcal{C}\left([-\rho,T];\mathcal{P}_2(\mathbb{R}^d)\right)
\end{equation*}
by $	\Phi_t(\mu)=\mathbb{L}_{\Upsilon^{\mu}(t)}$, where $\mathbb{L}_{\Upsilon^{\mu}(t)}$ is the distribution of $\Upsilon^{\mu}(t)$.
Next, we will show $\Phi$ is strictly contractive.
Using the same techniques of getting (3.5) in the proof of Theorem 3.1 in \cite{ccc12} gives that, for  $\bar{p}\geq2$, 
\begin{equation}\label{cccbbdt}
	\begin{split}
		&\mathbb{E}\big(\sup_{t\in[0,{T}]}\big|\int_{0}^{t}\beta_{s}(\mu_{s},\mu_{s-\rho})dB_s^H\big|^{\bar{p}}\big)\leq C_{\bar{p},T,H}\mathbb{E}\int_{0}^{{T}}\big|\beta_{t}(\mu_{t},\mu_{t-\rho})\big|^{\bar{p}}dt.
	\end{split}
\end{equation}
For any $t\in[0,T]$ and $ \hat{p}\in[2,\frac{\bar{p}}{2l}]$, we get from H\"older's inequality,  (\ref{cccbbdt}) and Assumption \ref{cccass1} that
\begin{equation*}
	\begin{split}
		&\mathbb{E}\big(\sup_{s\in[0,t]}|\Upsilon^{\mu}(s)-\Upsilon^{\nu}(s)|^{\hat{p}}\big)\\
		\leq&2^{\hat{p}-1}\mathbb{E}\Big(\sup_{s\in[0,t]}\Big|\int_{0}^{s}\big[\alpha_{r,\rho}^{\mu}(\Upsilon^{\mu}(r),\Upsilon^{\mu}(r-\rho))-\alpha_{r,\rho}^{\nu}(\Upsilon^{\nu}(r),\Upsilon^{\nu}(r-\rho))\big]dr\Big|^{\hat{p}}\Big)\\
		&+2^{\hat{p}-1}\mathbb{E}\Big(\sup_{s\in[0,t]}\big|\int_{0}^{s}\big[\beta_{r,\rho}^{\mu}-\beta_{r,\rho}^{\nu}\big]dB_r^H\Big|^{\hat{p}}\Big)\\
		\leq&C_{\hat{p},T,K_1}\mathbb{E}\int_{0}^{t}\Big[|\Upsilon^{\mu}(r)-\Upsilon^{\nu}(r)|^{\hat{p}}\\&+\big(1+|\Upsilon^{\mu}(r-\rho)|^l+|\Upsilon^{\nu}(r-\rho)|^l\big)^{\hat{p}}|\Upsilon^{\mu}(r-\rho)-\Upsilon^{\nu}(r-\rho)|^{\hat{p}}\\&+\mathbb{W}_{2}^{\hat{p}}(\mu_r,\nu_r)+\mathbb{W}_{2}^{\hat{p}}(\mu_{r-\rho},\nu_{r-\rho})\Big]dr
		+C_{\hat{p},T,H,K_2}\mathbb{E}\int_{0}^{t}\big|\beta_{r,\rho}^{\mu}-\beta_{r,\rho}^{\nu}\big|^{\hat{p}}dr\\
		\leq&c_{*}\int_{0}^{t}\Big[\mathbb{E}|\Upsilon^{\mu}(r)-\Upsilon^{\nu}(r)|^{\hat{p}}+\big[\mathbb{E}\big(1+|\Upsilon^{\mu}(r-\rho)|^{2\hat{p}l}+|\Upsilon^{\nu}(r-\rho)|^{2\hat{p}l}\big)\big]^\frac{1}{2}\big[\mathbb{E}|\Upsilon^{\mu}(r-\rho)-\Upsilon^{\nu}(r-\rho)|^{2\hat{p}}\big]^\frac{1}{2}\\&+\mathbb{W}_{2}^{\hat{p}}(\mu_r,\nu_r)+\mathbb{W}_{2}^{\hat{p}}(\mu_{r-\rho},\nu_{r-\rho})
		\Big]dr,
	\end{split}
\end{equation*}
where $c_{*}$ is a constant depending on $\hat{p},T,H,K_1,K_2$.
So, applying (\ref{cxcxa2}) with $\hat{p}\leq\frac{\bar{p}}{2l}$ means that 
\begin{equation}\label{cxcxa3}
	\begin{split}
		&\mathbb{E}\big(\sup_{s\in[0,t]}|\Upsilon^{\mu}(s)-\Upsilon^{\nu}(s)|^{\hat{p}}\big)\\
		\leq&c_{*}\int_{0}^{t}\mathbb{E}\big(\sup_{s\in[0,r]}|\Upsilon^{\mu}(s)-\Upsilon^{\nu}(s)|^{\hat{p}}\big)dr+c_{*}\int_{0}^{t}\big[\mathbb{E}\big(\sup_{s\in[0,r]}|\Upsilon^{\mu}(s-\rho)-\Upsilon^{\nu}(s-\rho)|^{2\hat{p}}\big)\big]^\frac{1}{2}dr\\&+c_{*}\int_{0}^{t}\sup_{s\in[0,r]}\mathbb{W}_{2}^{\hat{p}}(\mu_s,\nu_s)dr.
	\end{split}
\end{equation}
Define a sequence as
\begin{equation*}
	\hat{p}_{m}=(2-m+\lfloor\frac{T}{\rho}\rfloor)\hat{p}2^{1-m+\lfloor\frac{T}{\rho}\rfloor},~~~m=1, 2, \cdots, \lfloor\frac{T}{\rho}\rfloor+1.
\end{equation*}
Obviously, it holds that $\hat{p}_{\lfloor\frac{T}{\rho}\rfloor+1}=\hat{p}$ and $2\hat{p}_{m+1}<\hat{p}_{m}$ for $m=1, 2, \cdots, \lfloor\frac{T}{\rho}\rfloor$.

For $s\in[0,\rho]$, the Gronwall inequality with (\ref{cxcxa3}) gives that
\begin{equation}\label{cxcxa5}
	\begin{split}
		\sup_{s\in[0,\rho]}\mathbb{W}_{2}^{\hat{p}_1}(\Phi_s(\mu),\Phi_s(\nu))\leq
		\mathbb{E}\big(\sup_{s\in[0,\rho]}|\Upsilon^{\mu}(s)-\Upsilon^{\nu}(s)|^{\hat{p}_{1}}\big)
		\leq&c_{*}e^{c_{*}\rho}\int_{0}^{\rho}\sup_{s\in[0,r]}\mathbb{W}_{2}^{\hat{p}_1}(\mu_s,\nu_s)dr.
	\end{split}
\end{equation}
For $\rho_{*}>0$ (which is independent of the initial data) such that $c_{*}e^{c_{*}\rho_{*}}\rho_{*}\leq\frac{1}{2}$, let
\begin{equation*}
	\begin{split}
		\tilde{S}_{\rho_{*}}=\left\lbrace\mu_{.} \in\mathcal{C}\big([0,\rho_{*}];\mathcal{P}_2(\mathbb{R}^d)\big): \tilde{d}(\mu_{\cdot},\mu_{0})<\infty, \mu_{0}=\mathbb{L}_{X(0)} \right\rbrace 
	\end{split}
\end{equation*}
equipped with the uniform metric
\begin{equation*}
	\begin{split}
		\tilde{d}(\mu,\nu):=\sup_{s\in[0,\rho_{*}]}\mathbb{W}_{2}(\mu_s,\nu_s).
	\end{split}
\end{equation*}
We see that $(\tilde{S}_{\rho_{*}},\tilde{d})$ is a complete metric space.
Then using (\ref{cxcxa5}) implies
\begin{equation*}
	\begin{split}
			\tilde{d}^{\hat{p}_1}(\Phi(\mu),\Phi(\nu))\leq\frac{1}{2}	\tilde{d}^{\hat{p}_1}(\mu,\nu),
	\end{split}
\end{equation*}
which means that $\Phi$ is strictly contractive in the complete metric space $(\tilde{S}_{\rho_{*}},\tilde{d})$.
Hence, the Banach fixed point theorem with the definition of $\Phi$ reveals that there exists a unique $\mu\in\tilde{S}_{\rho_{*}}$ such that
\begin{equation*}
	\begin{split}
		\Phi_t(\mu)=\mu_t=\mathbb{L}_{\Upsilon^{\mu}(t)}
	\end{split}
\end{equation*}
on $t\in[0,\rho_{*}]$.
Thus, the strong wellposedness of (\ref{ccc1.1}) is obtained on $[0,\rho_{*}]$.
Repeating this procedure with the initial value $\Upsilon_{n\rho_{*}}$ (and the initial time $n\rho_{*}$) for $1\leq n \leq \lfloor\frac{\rho}{\rho_{*}}\rfloor+1$ gives the
strong wellposedness of (\ref{ccc1.1}) on $[0,\rho]$.

For $s\in[0,2\rho]$, the H\"older inequality with (\ref{cxcxa3}) and (\ref{cxcxa5}) gives that
\begin{equation}\label{cxcxa6}
	\begin{split}
		&\sup_{s\in[0,2\rho]}\mathbb{W}_{2}^{\hat{p}_2}(\Phi_s(\mu),\Phi_s(\nu))\\ \leq&
		\mathbb{E}\big(\sup_{s\in[0,2\rho]}|\Upsilon^{\mu}(s)-\Upsilon^{\nu}(s)|^{\hat{p}_{2}}\big)\\
		\leq&
		2c_{*}e^{2c_{*}\rho}\rho\left[\mathbb{E}\big(\sup_{s\in[0,\rho]}|\Upsilon^{\mu}(s)-\Upsilon^{\nu}(s)|^{\hat{p}_{1}}\big)\right]^{\frac{\hat{p}_2}{\hat{p}_1}}+2c_{*}e^{2c_{*}\rho}\rho\sup_{s\in[0,2\rho]}\mathbb{W}_{2}^{\hat{p}_2}(\mu_s,\nu_s)\\
		\leq&
		2c_{*}e^{2c_{*}\rho}\rho\left[c_{*}e^{c_{*}\rho}\rho \sup_{s\in[0,\rho]}\mathbb{W}_{2}^{\hat{p}_1}(\mu_s,\nu_s) \right]^{\frac{\hat{p}_2}{\hat{p}_1}}+2c_{*}e^{2c_{*}\rho}\rho\sup_{s\in[0,2\rho]}\mathbb{W}_{2}^{\hat{p}_2}(\mu_s,\nu_s)\\
		\leq&
		2c_{*}e^{2c_{*}\rho}\rho\left[1\vee(c_{*}e^{c_{*}\rho}\rho)^{\frac{\hat{p}_2}{\hat{p}_1}}\right]\sup_{s\in[0,2\rho]}\mathbb{W}_{2}^{\hat{p}_2}(\mu_s,\nu_s)\\
	\end{split}
\end{equation} 
For $\rho_{*}>0$ such that $2c_{*}e^{2c_{*}\rho_{*}}\rho_{*}[1\vee(c_{*}e^{c_{*}\rho_{*}}\rho_{*})^{\hat{p}_2/\hat{p}_1}]\leq\frac{1}{2}$, we get from (\ref{cxcxa6}) that
\begin{equation*}
	\begin{split}
			\tilde{d}^{\hat{p}_2}(\Phi(\mu),\Phi(\nu))\leq\frac{1}{2}\tilde{d}^{\hat{p}_2}(\mu,\nu),
	\end{split}
\end{equation*}
which means that $\Phi$ is strictly contractive in the complete metric space $(\tilde{S}_{\rho_{*}},\tilde{d})$.
Again, the Banach fixed point theorem gives that there exists a unique $\mu\in\tilde{S}_{\rho_{*}}$ such that
\begin{equation*}
	\begin{split}
		\Phi_t(\mu)=\mu_t=\mathbb{L}_{\Upsilon^{\mu}(t)}
	\end{split}
\end{equation*}
on $t\in[0,\rho_{*}]$.
Then the strong wellposedness of (\ref{ccc1.1}) is derived on $[0,\rho_{*}]$.
Repeating this procedure with the initial value $\Upsilon_{n\rho_{*}}$  for $1\leq n \leq \lfloor\frac{2\rho}{\rho_{*}}\rfloor+1$ gives the
strong wellposedness of (\ref{ccc1.1}) on $[0,2\rho]$.
Then by iteration about the time-delay segment, we can show the
strong wellposedness of (\ref{ccc1.1}) on $[0,T]$.

Then we give the detailed moment estimation  of the solution to (\ref{ccc1.1}).
For any $\bar{p}\geq 2$ and $t\in[0,T]$, the H\"older inequality leads to
\begin{equation*}
\begin{split}
	&\mathbb{E}\big(\sup_{r\in[0,t ]}|\Upsilon (r)|^{\bar{p}}\big)\\
	\leq& C_{\bar{p}}\mathbb{E}\|\xi\|^{\bar{p}}+C_{\bar{p},T}\mathbb{E}\int_{0}^{t }\big|\alpha_{r}(\Upsilon (r),\Upsilon (r-\rho),\mathbb{L}_{\Upsilon (r)},\mathbb{L}_{\Upsilon (r-\rho)})\big|^{\bar{p}}dr\\
	&+ C_{\bar{p}}\mathbb{E}\big(\sup_{r\in[0,t ]}\big|\int_{0}^{r}\beta_{u}(\mathbb{L}_{\Upsilon (u)},\mathbb{L}_{\Upsilon (u-\rho)})dB_u^H\big|^{\bar{p}}\big).
\end{split}
\end{equation*}

We derive from Assumption \ref{cccass1}, Lemma \ref{ccclmm1} and (\ref{cccbbdt}) that
\begin{equation*}
	\begin{split}
		&\mathbb{E}\big(\sup_{r\in[0,t ]}|\Upsilon (r)|^{\bar{p}}\big)\\
		\leq& C_{\bar{p}}\mathbb{E}\|\xi\|^{\bar{p}}+C_{\bar{p},T,K_1}\mathbb{E}\int_{0}^{t }\big[1+|\Upsilon (r)|^{\bar{p}}+|\Upsilon (r-\rho)|^{\bar{p}(l+1)}+\mathcal{W}_{2}^{\bar{p}}(\mathbb{L}_{\Upsilon (r)})+\mathcal{W}_{2}^{\bar{p}}(\mathbb{L}_{\Upsilon (r-\rho)})\big]dr\\
		&+ C_{\bar{p},T,H,K_2}\mathbb{E}\int_{0}^{t }\big[1+\mathcal{W}_{2}^{\bar{p}}(\mathbb{L}_{\Upsilon (r)})+\mathcal{W}_{2}^{\bar{p}}(\mathbb{L}_{\Upsilon (r-\rho)})\big]dr\\
			\leq& C_{\bar{p},T,H,K_1,K_2,\|\xi\|}\big[1+\int_{0}^{t}\mathbb{E}\big(\sup_{u\in[0,r ]}|\Upsilon (u)|^{\bar{p}}\big)dr+\mathbb{E}\int_{0}^{t }|\Upsilon (r-\rho)|^{\bar{p}(l+1)} dr\big].
	\end{split}
\end{equation*}
Thanks to the Gronwall inequality,
\begin{equation*}
	\begin{split}
		&\mathbb{E}\big(\sup_{t\in[0,T ]}|\Upsilon (t)|^{\bar{p}}\big)
		\leq C_{\bar{p},T,H,K_1,K_2,\|\xi\|}\big[1+\mathbb{E}\big(\sup_{t\in[0,T ]}|\Upsilon (t-\rho)|^{\bar{p}(l+1)}\big)\big].
	\end{split}
\end{equation*}
Constructing a sequence $	\bar{p}_{m}=(2-m+\lfloor\frac{T}{\rho}\rfloor)\bar{p}(l+1)^{1-m+\lfloor\frac{T}{\rho}\rfloor}$ and using the same iteration technique about $\bar{p}_{m}$ as that in Lemma \ref{cccsdeexistence} yield that
\begin{equation*}
	\begin{split}
		&\mathbb{E}\big(\sup_{t\in[0,((\lfloor\frac{T}{\rho}\rfloor+1)\rho) ]}|\Upsilon (t)|^{\bar{p}}\big)
		\leq C_{\bar{p},T,H,K_1,K_2,\|\xi\|,\bar{p}_{m}},
	\end{split}
\end{equation*}
where $m=1, 2, \cdots, \lfloor\frac{T}{\rho}\rfloor+1$.
So the desired result holds.

For the case $H\in (0,1/2)$, since $\beta$ does not depend on the distribution, the assertion can be obtained analogically by means of  Lemma \ref{cccineq1}.
\eproof

In the rest of this paper, the DSMVEs  are all autonomous for convenience.
Moreover, 
the relevant settings of the equations are stated for the case  $H\in (1/2,1)$. As for $H\in (0,1/2)$, we just let the coefficient $\beta(\cdot,\cdot)$ become a constant $\beta$.

\section{Propagation of chaos}\label{cccsec03}

In this section, the stochastic particle method in \cite{ccc25,ccc26} is used to approximate the distribution in the coefficient of DSMVE.
For any $i\in \mathbb{S}_N$, consider the non-interacting particle system
\begin{equation} \label{cccnonIPS}
	d\Upsilon^{i}(t)
	=\alpha\left(\Upsilon^{i}(t),\Upsilon^{i}(t-\rho),\mathbb{L}_{\Upsilon^{i}(t)},\mathbb{L}_{\Upsilon^{i}(t-\rho)}\right)dt+\beta\left(\mathbb{L}_{\Upsilon^{i}(t)},\mathbb{L}_{\Upsilon^{i}(t-\rho)}\right) dB_t^{H,i},~~~ t\in[0,T],
\end{equation}
with the initial value $X_{0}^{i}=\xi^{i}$, which is an $\mathcal{F}_0$-measurable $\mathcal{C}$-valued
random variable with $\mathbb{E}\|\xi\|^{\check{p}}<\infty$ for any $\check{p}>0$, where $\mathbb{L}_{\Upsilon^{i}(\cdot)}$ is distribution of $\Upsilon^{i}(\cdot)$.
Here, $(\xi^{i},B^{H,i})$ are the independent copies of $(\xi,B^H)$, and all $(\xi^{i},B^{H,i})$ are  independent and identically distributed. Moreover, it holds that $\mathbb{L}_{\Upsilon^{i}(t)}=\mathbb{L}_{\Upsilon(t)}$ for any $t\in [0,T]$ and $i\in\mathbb{S}_N$.
The corresponding interacting particle system driven by fBm is
\begin{equation} \label{cccIPS}
	d\Upsilon^{i,N}(t)
	=\alpha\left(\Upsilon^{i,N}(t),\Upsilon^{i,N}(t-\rho),\mathbb{L}_{\Upsilon^{N}(t)},\mathbb{L}_{\Upsilon^{N}(t-\rho)}\right)dt+\beta\left(\mathbb{L}_{\Upsilon^{N}(t)},\mathbb{L}_{\Upsilon^{N}(t-\rho)}\right) dB_t^{H,i},~~~ t\in[0,T],
\end{equation}
with the initial value $X_{0}^{i}=\xi^{i}$, where $\mathbb{L}_{\Upsilon^{N}(\cdot)}:=\frac{1}{N}\sum_{j=1}^{N}\delta_{\Upsilon^{j,N}(\cdot)}$.
Under Assumption \ref{cccass1}, one can easily get the wellposedness of (\ref{cccnonIPS}) and (\ref{cccIPS}).

\begin{lemma}
For $H\in (0,1/2)\cup (1/2,1)$, let Assumption \ref{cccass1} hold. Then for any $\bar{p}\geq2$ and $T>0$,
	\begin{equation*}\label{cccIPSbound}
		\mathbb{E}\big(\sup_{t\in[0,T]}|\Upsilon^{i}(t)|^{\bar{p}}\big)\vee	\mathbb{E}\big(\sup_{t\in[0,T]}|\Upsilon^{i,N}(t)|^{\bar{p}}\big)\leq C_{\bar{p},T,H,K_1,K_2,\|\xi\|,l,\bar{p}_m}.
	\end{equation*}
\end{lemma}

\begin{theorem}\label{cccfifi1}
	When $H\in (0,1/2)\cup (1/2,1)$, let Assumption \ref{cccass1} hold with $p(pl+\varepsilon)<\varepsilon\bar{p}$ for $\varepsilon\in(0,1]$. Then we derive that, for any $i\in \mathbb{S}_N$ and $p\geq 2$,
	\begin{equation*}
		\begin{split}
			\mathbb{E}\Big(\sup_{t\in[0,T]}|\Upsilon^{i}(t)-\Upsilon^{i,N}(t)|^{p}\Big)\leq C
			\left\{\begin{array}{ll}
				(N^{-1 / 2})^{\lambda_{p,T,\rho}}, & \text { if } p>d/2, \\
				{[N^{-1 / 2} \log(1+N)]}^{\lambda_{p,T,\rho}}, & \text { if } p=d/2,\\
				{(N^{-p / d})}^{\lambda_{p,T,\rho}}, & \text { if }2\leq p<d/2,
			\end{array}\right.
		\end{split}
	\end{equation*}
	where $\lambda_{p,T,\rho}=(\frac{p-\varepsilon}{p})^{\lfloor\frac{T}{\rho}\rfloor}$ and $C$ is a positive  constant dependent of $\bar{p}$, $T$, $H$, $K_1$, $K_2$, $\|\xi\|$, $ \bar{p}_m$ but independent of $N$.
\end{theorem}

\noindent
{\it Proof}.
We begin with the case $H\in (1/2,1)$.
For any $i\in\mathbb{S}_N$, $p\geq 2$ and $t\in[0,T]$, we derive from H\"{o}lder's inequality and (\ref{cccbbdt}) that
\begin{equation*}
	\begin{split}
		&\mathbb{E}\big(\sup_{s\in[0,t]}|\Upsilon^{i}(s)-\Upsilon^{i,N}(s)|^p\big)\\
		\leq&2^{p-1}\mathbb{E}\Big(\sup_{s\in[0,t]}\Big|\int_{0}^{s}\big[\alpha\big(\Upsilon^{i}(r),\Upsilon^{i}(r-\rho),\mathbb{L}_{\Upsilon^{i}(r)},\mathbb{L}_{\Upsilon^{i}(r-\rho)}\big)\\&-\alpha\big(\Upsilon^{i,N}(r),\Upsilon^{i,N}(r-\rho),\mathbb{L}_{\Upsilon^{N}(r)},\mathbb{L}_{\Upsilon^{N}(r-\rho)}\big)\big]dr\Big|^p\Big)\\
		&+2^{p-1}\mathbb{E}\Big(\sup_{s\in[0,t]}\Big|\int_{0}^{s}\big[\beta\big(\mathbb{L}_{\Upsilon^{i}(r)},\mathbb{L}_{\Upsilon^{i}(r-\rho)}\big)-\beta\big(\mathbb{L}_{\Upsilon^{N}(r)},\mathbb{L}_{\Upsilon^{N}(r-\rho)}\big)\big]dB_r^H\Big|^{\hat{p}}\Big)\\
		\leq&C_{p,T,K_1}\mathbb{E}\int_{0}^{t}\Big[|\Upsilon^{i}(r)-\Upsilon^{i,N}(r)|\\&+\big(1+|\Upsilon^{i}(r-\rho)|^l+|\Upsilon^{i,N}(r-\rho)|^l\big)|\Upsilon^{i}(r-\rho)-\Upsilon^{i,N}(r-\rho)|\\&+\mathbb{W}_{2}(\mathbb{L}_{\Upsilon^{i}(r)},\mathbb{L}_{\Upsilon^{N}(r)})+\mathbb{W}_{2}(\mathbb{L}_{\Upsilon^{i}(r-\rho)},\mathbb{L}_{\Upsilon^{N}(r-\rho)})\Big]^pdr\\
		&+C_{p,T,H,K_2}\mathbb{E}\int_{0}^{t}\Big[\mathbb{W}_{2}(\mathbb{L}_{\Upsilon^{i}(r)},\mathbb{L}_{\Upsilon^{N}(r)})+\mathbb{W}_{2}(\mathbb{L}_{\Upsilon^{i}(r-\rho)},\mathbb{L}_{\Upsilon^{N}(r-\rho)})\Big]^pdr\\
		\leq&C_{p,T,H,K_1,K_2}\mathbb{E}\int_{0}^{t}\Big[|\Upsilon^{i}(r)-\Upsilon^{i,N}(r)|^p\\&+\big(1+|\Upsilon^{i}(r-\rho)|^l+|\Upsilon^{i,N}(r-\rho)|^l\big)^p|\Upsilon^{i}(r-\rho)-\Upsilon^{i,N}(r-\rho)|^p\\&+
	\mathbb{W}_{p}^{p}(\mathbb{L}_{\Upsilon^{i}(r)},\mathbb{L}_{\Upsilon^{N}(r)})+\mathbb{W}_{p}^{p}(\mathbb{L}_{\Upsilon^{i}(r-\rho)},\mathbb{L}_{\Upsilon^{N}(r-\rho)})	\Big]dr.
	\end{split}
\end{equation*}
Due to H\"older's inequality and Lemma \ref{cccIPSbound}, for $\frac{p(lp+\varepsilon)}{\varepsilon}<\bar{p}$ with $\varepsilon\in(0,1]$, we see that
\begin{equation*}
	\begin{split}
		&\mathbb{E}\big[\big(1+|\Upsilon^{i}(r-\rho)|^l+|\Upsilon^{i,N}(r-\rho)|^l\big)^p|\Upsilon^{i}(r-\rho)-\Upsilon^{i,N}(r-\rho)|^p\big]\\
		\leq&C_p\big[\mathbb{E}\big(1+|\Upsilon^{i}(r-\rho)|^{lp+\varepsilon}+|\Upsilon^{i,N}(r-\rho)|^{lp+\varepsilon}\big)^\frac{p}{\varepsilon}\big]^\frac{\varepsilon}{p}\big[\mathbb{E}|\Upsilon^{i}(r-\rho)-\Upsilon^{i,N}(r-\rho)|^p\big]^\frac{p-\varepsilon}{p}\\
		\leq&C_{p,T,H,K_1,K_2,\|\xi\|,\bar{p}_m}[\mathbb{E}|\Upsilon^{i}(r-\rho)-\Upsilon^{i,N}(r-\rho)|^p\big]^\frac{p-\varepsilon}{p}.
	\end{split}
\end{equation*}
The Gronwall inequality leads to
\begin{equation}\label{cccpoc1}
	\begin{split}
		&\mathbb{E}\big(\sup_{s\in[0,t]}|\Upsilon^{i}(s)-\Upsilon^{i,N}(s)|^p\big)\\
		\leq&C_{p,T,H,K_1,K_2,\|\xi\|,\bar{p}_m}\Big(\big[\mathbb{E}\big(\sup_{s\in[0,t]}|\Upsilon^{i}(s-\rho)-\Upsilon^{i,N}(s-\rho)|^p\big)\big]^\frac{p-\varepsilon}{p}\\&+\mathbb{E}\int_{0}^{t}\big(
		\mathbb{W}_{p}^{p}(\mathbb{L}_{\Upsilon^{i}(r)},\mathbb{L}_{\Upsilon^{N}(r)})+\mathbb{W}_{p}^{p}(\mathbb{L}_{\Upsilon^{i}(r-\rho)},\mathbb{L}_{\Upsilon^{N}(r-\rho)})\big)dr	\Big).
	\end{split}
\end{equation}
For $s\in[0,\rho]$, we know from (\ref{cccpoc1}) that
\begin{equation*}
	\begin{split}
		\mathbb{E}\big(\sup_{s\in[0,\rho]}|\Upsilon^{i}(s)-\Upsilon^{i,N}(s)|^p\big)
		\leq C_{p,T,H,K_1,K_2,\|\xi\|,\bar{p}_m}\mathbb{E}\int_{0}^{\rho}\big(
		\mathbb{W}_{p}^{p}(\mathbb{L}_{\Upsilon^{i}(r)},\mathbb{L}_{\Upsilon^{N}(r)})+\mathbb{W}_{p}^{p}(\mathbb{L}_{\Upsilon^{i}(r-\rho)},\mathbb{L}_{\Upsilon^{N}(r-\rho)})\big)dr.
	\end{split}
\end{equation*}
Define the empirical measure $\hat{\mathbb{L}}_{\Upsilon^{i,N}(\cdot)}=\frac{1}{N}\sum_{j=1}^{N}\delta_{\Upsilon^{i}(\cdot)}$.
For any $i\in\mathbb{S}_N$, 
\begin{equation*}
	\begin{split}
	\mathbb{W}_{p}^{p}(\mathbb{L}_{\Upsilon^{i}(\cdot)},\mathbb{L}_{\Upsilon^{N}(\cdot)})
		\leq &C_p
		\mathbb{W}_{p}^{p}(\mathbb{L}_{\Upsilon^{i}(\cdot)},\hat{\mathbb{L}}_{\Upsilon^{i,N}(\cdot)})+C_p\mathbb{W}_{p}^{p}(\hat{\mathbb{L}}_{\Upsilon^{i,N}(\cdot)},\mathbb{L}_{\Upsilon^{N}(\cdot)})\\
			\leq &C_p
		\mathbb{W}_{p}^{p}(\mathbb{L}_{\Upsilon^{i}(\cdot)},\hat{\mathbb{L}}_{\Upsilon^{i,N}(\cdot)})+C_p\frac{1}{N}\sum_{j=1}^{N}|\Upsilon^{j}(\cdot)-\Upsilon^{j,N}(\cdot)|^p.
	\end{split}
\end{equation*}
Since the distributions of all $j$ are identical, we obtain
\begin{equation*}
	\begin{split}
		&\mathbb{E}\big(\sup_{s\in[0,\rho]}|\Upsilon^{i}(s)-\Upsilon^{i,N}(s)|^p\big)\\
		\leq&C_{p,T,H,K_1,K_2,\|\xi\|,\bar{p}_m}\mathbb{E}\int_{0}^{\rho}\Big[	\mathbb{W}_{p}^{p}(\mathbb{L}_{\Upsilon^{i}(r)},\hat{\mathbb{L}}_{\Upsilon^{i,N}(r)})+\mathbb{W}_{p}^{p}(\mathbb{L}_{\Upsilon^{i}(r-\rho)},\hat{\mathbb{L}}_{\Upsilon^{i,N}(r-\rho)})\\&+|\Upsilon^{i}(r)-\Upsilon^{i,N}(r)|^p+|\Upsilon^{i}(r-\rho)-\Upsilon^{i,N}(r-\rho)|^p\Big]dr	\\
			\leq&C_{p,T,H,K_1,K_2,\|\xi\|,\bar{p}_m}\int_{0}^{\rho}\mathbb{E}\big(\sup_{s\in[0,r]}|\Upsilon^{i}(s)-\Upsilon^{i,N}(s)|^p\big)dr\\
			&+C_{p,T,H,K_1,K_2,\|\xi\|,\bar{p}_m}\mathbb{E}\int_{0}^{\rho}\big[	\mathbb{W}_{p}^{p}(\mathbb{L}_{\Upsilon^{i}(r)},\hat{\mathbb{L}}_{\Upsilon^{i,N}(r)})+\mathbb{W}_{p}^{p}(\mathbb{L}_{\Upsilon^{i}(r-\rho)},\hat{\mathbb{L}}_{\Upsilon^{i,N}(r-\rho)})\big]dr.
	\end{split}
\end{equation*}
In view of Theorem 1 in \cite{ccc29} with Gronwall's inequality, we get
$$\mathbb{E}\big(\sup_{s\in[0,\rho]}|\Upsilon^{i}(s)-\Upsilon^{i,N}(s)|^p\big)\leq C_{d,p,T,H,K_1,K_2,\|\xi\|,\bar{p}_m}\begin{cases}N^{-\frac{1}{2}}+N^{-\frac{\bar{p}-p}{\bar{p}}},& \text{if}~~ p>\frac{d}{2},\bar{p}\neq 2p,\\N^{-\frac{1}{2}}\log(1+N)+N^{-\frac{\bar{p}-p}{\bar{p}}},&\text{if}~~ p=\frac{d}{2},\bar{p}\neq 2p,\\N^{-\frac{p}{d}}+N^{-\frac{\bar{p}-p}{\bar{p}}},&\text{if}~~ 2\leq p<\frac{d}{2}.\end{cases}$$
Actually, the condition $\varepsilon\bar{p}>p(lp+\varepsilon)$ with $p\geq2$ means $\bar{p}>2p\geq4$. So the above inequality becomes
\begin{equation}\label{cccbigpoc}
	\begin{split}
		\mathbb{E}\big(\sup_{s\in[0,\rho]}|\Upsilon^{i}(s)-\Upsilon^{i,N}(s)|^p\big)\leq C_{d,p,T,H,K_1,K_2,\|\xi\|,\bar{p}_m}\begin{cases}N^{-\frac{1}{2}},& \text{if}~~ p>\frac{d}{2},\\N^{-\frac{1}{2}}\log(1+N),&\text{if}~~ p=\frac{d}{2},\\N^{-\frac{p}{d}},&\text{if}~~ 2\leq p<\frac{d}{2}.\end{cases}	
	\end{split}
	\end{equation}
For $s\in[0,2\rho]$, the H\"older inequality with (\ref{cccpoc1}) leads to
\begin{equation*}
	\begin{split}
		&\mathbb{E}\big(\sup_{s\in[0,2\rho]}|\Upsilon^{i}(s)-\Upsilon^{i,N}(s)|^p\big)\\
		\leq&C_{p,T,H,K_1,K_2,\|\xi\|,\bar{p}_m}\big[\mathbb{E}\big(\sup_{s\in[0,2\rho]}|\Upsilon^{i}(s-\rho)-\Upsilon^{i,N}(s-\rho)|^p\big)\big]^\frac{p-\varepsilon}{p}\\&+C_{p,T,H,K_1,K_2,\|\xi\|,\bar{p}_m}\mathbb{E}\int_{0}^{2\rho}
	\Big[	\mathbb{W}_{p}^{p}(\mathbb{L}_{\Upsilon^{i}(r)},\hat{\mathbb{L}}_{\Upsilon^{i,N}(r)})+\mathbb{W}_{p}^{p}(\hat{\mathbb{L}}_{\Upsilon^{i,N}(r)},\mathbb{L}_{\Upsilon^{N}(r)})\\&+\mathbb{W}_{p}^{p}(\mathbb{L}_{\Upsilon^{i}(r-\rho)},\hat{\mathbb{L}}_{\Upsilon^{i,N}(r-\rho)})+\mathbb{W}_{p}^{p}(\hat{\mathbb{L}}_{\Upsilon^{i,N}(r-\rho)},\mathbb{L}_{\Upsilon^{N}(r-\rho)})	\Big]dr\\
		\leq&C_{p,T,H,K_1,K_2,\|\xi\|,\bar{p}_m}\big[\mathbb{E}\big(\sup_{s\in[0,\rho]}|\Upsilon^{i}(s)-\Upsilon^{i,N}(s)|^p\big)\big]^\frac{p-\varepsilon}{p}\\&+C_{p,T,H,K_1,K_2,\|\xi\|,\bar{p}_m}\mathbb{E}\int_{0}^{2\rho}
	\Big[	\mathbb{W}_{p}^{p}(\mathbb{L}_{\Upsilon^{i}(r)},\hat{\mathbb{L}}_{\Upsilon^{i,N}(r)})+\mathbb{W}_{p}^{p}(\mathbb{L}_{\Upsilon^{i}(r-\rho)},\hat{\mathbb{L}}_{\Upsilon^{i,N}(r-\rho)})\Big]dr\\
	&+C_{p,T,H,K_1,K_2,\|\xi\|,\bar{p}_m}\int_{0}^{2\rho}\mathbb{E}\big(\sup_{s\in[0,r]}|\Upsilon^{i}(s)-\Upsilon^{i,N}(s)|^p\big)dr.
	\end{split}
\end{equation*}
Using Gronwall's inequality with (\ref{cccbigpoc}) on $s\in[0,\rho]$ gives that
\begin{equation*}
	\begin{split}
		&\mathbb{E}\big(\sup_{s\in[0,2\rho]}|\Upsilon^{i}(s)-\Upsilon^{i,N}(s)|^p\big)\\
		\leq&C_{p,T,H,K_1,K_2,\|\xi\|,\bar{p}_m}\big[\mathbb{E}\big(\sup_{s\in[0,\rho]}|\Upsilon^{i}(s)-\Upsilon^{i,N}(s)|^p\big)\big]^\frac{p-\varepsilon}{p}\\&+C_{p,T,H,K_1,K_2,\|\xi\|,\bar{p}_m}\mathbb{E}\int_{0}^{2\rho}
		\Big[	\mathbb{W}_{p}^{p}(\mathbb{L}_{\Upsilon^{i}(r)},\hat{\mathbb{L}}_{\Upsilon^{i,N}(r)})+\mathbb{W}_{p}^{p}(\mathbb{L}_{\Upsilon^{i}(r-\rho)},\hat{\mathbb{L}}_{\Upsilon^{i,N}(r-\rho)})\Big]dr\\
	\leq &C_{d,p,T,H,K_1,K_2,\|\xi\|,\bar{p}_m}\begin{cases}(N^{-\frac{1}{2}})^\frac{p-\varepsilon}{p},& \text{if}~~ p>\frac{d}{2},\\
		\big(N^{-\frac{1}{2}}\log(1+N)\big)^\frac{p-\varepsilon}{p},&\text{if}~~ p=\frac{d}{2},\\
		(N^{-\frac{p}{d}})^\frac{p-\varepsilon}{p},&\text{if}~~ 2\leq p<\frac{d}{2}.\end{cases}	
	\end{split}
\end{equation*}
For $s\in[0,3\rho]$, it is similar to see
\begin{equation*}
	\begin{split}
		&\mathbb{E}\big(\sup_{s\in[0,3\rho]}|\Upsilon^{i}(s)-\Upsilon^{i,N}(s)|^p\big)\\
		\leq&C_{p,T,H,K_1,K_2,\|\xi\|,\bar{p}_m}\big[\mathbb{E}\big(\sup_{s\in[0,2\rho]}|\Upsilon^{i}(s)-\Upsilon^{i,N}(s)|^p\big)\big]^\frac{p-\varepsilon}{p}\\&+C_{p,T,H,K_1,K_2,\|\xi\|,\bar{p}_m}\mathbb{E}\int_{0}^{3\rho}
		\Big[	\mathbb{W}_{p}^{p}(\mathbb{L}_{\Upsilon^{i}(r)},\hat{\mathbb{L}}_{\Upsilon^{i,N}(r)})+\mathbb{W}_{p}^{p}(\mathbb{L}_{\Upsilon^{i}(r-\rho)},\hat{\mathbb{L}}_{\Upsilon^{i,N}(r-\rho)})\Big]dr\\
		\leq &C_{d,p,T,H,K_1,K_2,\|\xi\|,\bar{p}_m}\begin{cases}(N^{-\frac{1}{2}})^{(\frac{p-\varepsilon}{p})^2},& \text{if}~~ p>\frac{d}{2},\\
			\big(N^{-\frac{1}{2}}\log(1+N)\big)^{(\frac{p-\varepsilon}{p})^2},&\text{if}~~ p=\frac{d}{2},\\
			(N^{-\frac{p}{d}})^{(\frac{p-\varepsilon}{p})^2},&\text{if}~~ 2\leq p<\frac{d}{2}.\end{cases}	
	\end{split}
\end{equation*}
The desired result follows by the iteration about the time segment generated by the delay $\rho$.

As for the case $H\in (0,1/2)$, the stochastic integral in $\Upsilon^{i}(s)-\Upsilon^{i,N}(s)$ vanishes since it is an additive noise, then the target can be achieved through the similar process.
\eproof

\section{Numerical scheme}\label{cccsec04}

In this section, the classical EM scheme, which is not modified, is established for interacting particle system (\ref{cccIPS}) whose delay term is superlinear.
Let $\Delta=\frac{\rho}{M}=\frac{T}{M_T}$ for some positive integers $M$ and $M_T$. Set $t_k=k\Delta$ for $k=-M,\cdots,0,\cdots,M_T$ and define the EM scheme as
\begin{equation}\label{cccemli}
		\begin{split}
Z ^{i,N}(t_{k+1})&=Z ^{i,N}(t_k)+\alpha\big(Z ^{i,N}(t_k),Z ^{i,N}(t_{k-M}),\mathbb{L}_{Z ^{N}(t_k)},\mathbb{L}_{Z ^{N}(t_{k-M})}\big)\Delta\\&+\beta\big(\mathbb{L}_{Z ^{N}(t_k)},\mathbb{L}_{Z ^{N}(t_{k-M})}\big)\Delta B_{k}^{H,i},~~~k=0,1,\cdots,M_T-1,
	\end{split}
\end{equation}
with the initial value $Z ^{i,N}(t_k)=\xi(t_k)$, $k=-M,\cdots,1,0$, where $\Delta B_{k}^{H,i}=B_{t_{k+1}}^{H,i}-B_{t_{k}}^{H,i}$ and $ \mathbb{L}_{Z ^{N}(\cdot)}=\frac{1}{N}\sum_{j=1}^{N}\delta_{Z ^{j,N}(\cdot)}$.

For $t\in[0,T]$, the continuous-sample numerical scheme is 
\begin{equation}\label{ccccon1}
	\begin{split}
		\tilde{Z }^{i,N}(t)&=\xi^{i}(0)+\int_{0}^{t}\alpha\big(Z ^{i,N}(s),Z ^{i,N}(s-\rho),\mathbb{L}_{Z ^{N}(s)},\mathbb{L}_{Z ^{N}(s-\rho)}\big)ds\\&+\int_{0}^{t}\beta\big(\mathbb{L}_{Z ^{N}(s)},\mathbb{L}_{Z ^{N}(s-\rho)}\big)dB_{s}^{H,i},
	\end{split}
\end{equation}
with the step process $Z ^{i,N}(t):=\sum_{k=-M}^{M_T}Z ^{i,N}(t_k)\mathbb{I}_{[t_k,t_{k+1})}(t)$ and the empirical measure $ \mathbb{L}_{Z ^{N}(t)}=\frac{1}{N}\sum_{j=1}^{N}\delta_{Z ^{j,N}(t)}$.
Obviously,  $\tilde{Z }^{i,N}(t_k)=Z ^{i,N}(t_k)=Z ^{i,N}(t)$ for $t\in[t_k,t_{k+1})$.
And the diffusion coefficient in numerical scheme becomes a constant $\beta$ for $H\in (0,1/2)$.
\begin{lemma}\label{cccnumbound}
	For $H\in (0,1/2)\cup (1/2,1)$, let Assumption \ref{cccass1} hold. Then, for $\bar{p}\geq2$,
	\begin{equation*}
		\begin{split}
		\max_{i\in\mathbb{S}_N} \sup_{\Delta\in(0,1]} \mathbb{E}\big(\sup_{t\in[0,T]}|\tilde{Z }^{i,N}(t)|^{\bar{p}}\big) \leq C_{\bar{p},T,H,K_1,K_2,\|\xi\|,\bar{p}_m}.
		\end{split}
	\end{equation*} 
\end{lemma}

\noindent
{\it Proof}. 
We first discuss the case $H\in (1/2,1)$.
For any $t\in[0,T]$ and $i\in\mathbb{S}_N$, the H\"older inequality with Assumption \ref{cccass1} and (\ref{cccbbdt}) yields that
\begin{equation*}
	\begin{split}
	&\mathbb{E}\big(\sup_{s\in[0,t]}|\tilde{Z }^{i,N}(s)|^{\bar{p}}\big)\\\leq&3^{\bar{p}-1}\mathbb{E}\|\xi^{i}(0)\|^{\bar{p}}+3^{\bar{p}-1}\mathbb{E}\Big(\sup_{s\in[0,t]}\big|\int_{0}^{s}\alpha\big(Z ^{i,N}(r),Z ^{i,N}(r-\rho),\mathbb{L}_{Z ^{N}(r)},\mathbb{L}_{Z ^{N}(r-\rho)}\big)dr\big|^{\bar{p}}\Big)\\&+3^{\bar{p}-1}\mathbb{E}\Big(\sup_{s\in[0,t]}\big|\int_{0}^{s}\beta\big(\mathbb{L}_{Z ^{N}(r)},\mathbb{L}_{Z ^{N}(r-\rho)}\big)dB_{r}^{H,i}\big|^{\bar{p}}\Big)\\
	\leq& 3^{\bar{p}-1}\mathbb{E}\|\xi\|^{\bar{p}}+C_{\bar{p},T,H,K_1,K_2}\mathbb{E}\int_{0}^{t}\big[1+|Z ^{i,N}(r)|^{\bar{p}}+|Z ^{i,N}(r-\rho)|^{\bar{p}(l+1)}\\&+\mathcal{W}_{2}^{\bar{p}}(\mathbb{L}_{Z ^{N}(r)})+\mathcal{W}_{2}^{\bar{p}}(\mathbb{L}_{Z ^{N}(r-\rho)})\big]dr\\
	\leq& C_{\bar{p},T,H,K_1,K_2,\|\xi\|}\int_{0}^{t}\big(1+\mathbb{E}\big(\sup_{r\in[0,s]}|\tilde{Z }^{i,N}(r)|^{\bar{p}}\big)\big)ds+C_{\bar{p},T,H,K_1,K_2}\mathbb{E}\big(\sup_{s\in[0,t]}|\tilde{Z }^{i,N}(s-\rho)|^{\bar{p}(l+1)}\big),
	\end{split}
\end{equation*}
Thanks to Gronwall's inequality, we see 
\begin{equation*}
	\begin{split}
	\mathbb{E}\big(\sup_{s\in[0,t]}|\tilde{Z }^{i,N}(s)|^{\bar{p}}\big)
		\leq C_{\bar{p},T,H,K_1,K_2,\|\xi\|}\big[1+\mathbb{E}\big(\sup_{s\in[0,t]}|\tilde{Z }^{i,N}(s-\rho)|^{\bar{p}(l+1)}\big)\big].
	\end{split}
\end{equation*}
Construct a sequence: 
$
	\bar{p}_{m}=(2-m+\lfloor\frac{T}{\rho}\rfloor)\bar{p}(l+1)^{1-m+\lfloor\frac{T}{\rho}\rfloor}
$ for $m=1, 2, \cdots, \lfloor\frac{T}{\rho}\rfloor+1$.
For any $i\in\mathbb{S}_N$, then using the same iteration technique  as that in Lemma \ref{cccsdeexistence} leads to 
\begin{equation*}
	\begin{split}
		\mathbb{E}\big(\sup_{s\in[0,t]}|\tilde{Z }^{i,N}(s)|^{\bar{p}}\big)
		\leq C_{\bar{p},T,H,K_1,K_2,\|\xi\|,\bar{p}_m}.
	\end{split}
\end{equation*}

As for the case $H\in (0,1/2)$,
the proof process is similar after using  Lemma \ref{cccineq1}.

\eproof

\begin{lemma}\label{ccccondis}
	For $H\in  (1/2,1)$, let Assumption \ref{cccass1} hold. For any $i\in\mathbb{S}_N$ and $\hat{p}\in[2,\frac{\bar{p}}{l+1}]$ with $\bar{p}>2l+2$, we derive
	\begin{equation*}
		\begin{split}
			 \mathbb{E}\big(\sup_{t\in[0,T]}|\tilde{Z }^{i,N}(t)-Z ^{i,N}(t)|^{\hat{p}}\big) \leq C_{\hat{p},\bar{p},T,H,K_1,K_2,\|\xi\|,\bar{p}_m}\Delta^{\hat{p}H}.
		\end{split}
	\end{equation*} 
\end{lemma}

\noindent
{\it Proof}.
  For any $t\in[t_k,t_{k+1})$ and $i\in\mathbb{S}_N$, we obtain from (\ref{ccccon1}) that
  \begin{equation*}
  	\begin{split}
  		&\tilde{Z }^{i,N}(t)-Z ^{i,N}(t)\\=&\int_{t_k}^{t}\alpha\big(Z ^{i,N}(s),Z ^{i,N}(s-\rho),\mathbb{L}_{Z ^{N}(s)},\mathbb{L}_{Z ^{N}(s-\rho)}\big)ds+\int_{t_k}^{t}\beta\big(\mathbb{L}_{Z ^{N}(s)},\mathbb{L}_{Z ^{N}(s-\rho)}\big)dB_{s}^{H,i}\\
  		=:&J_1(t)+J_2(t).
  	\end{split}
  \end{equation*}
We first estimate $J_2(t)$. For $\hat{p}\geq2$, choose $\eta$ to satisfy $1-H<\eta<1-\frac{1}{\hat{p}}$.
Denote $\varphi_{\eta}=\int_{s}^{t}(t-u)^{-\eta}(u-s)^{\eta-1}du$.
Applying stochastic Fubini's theorem and H\"older's inequality yields that
  \begin{equation*}
  	\begin{split}
  		&\mathbb{E}\big(\sup_{t\in[t_k,t_{k+1}]}|J_2(t)|^{\hat{p}}\big)\\=&{(\varphi_{\eta})}^{-\hat{p}}\mathbb{E}\Big(\sup_{t\in[t_k,t_{k+1}]}\big|\int_{t_k}^{t} \big(\int_{s}^{t}(t-u)^{-\eta}(u-s)^{\eta-1}du\big) \beta\big(\mathbb{L}_{Z ^{N}(s)},\mathbb{L}_{Z ^{N}(s-\rho)}\big)        dB_{s}^{H,i}\big|^{\hat{p}}\Big)\\
  		=&{(\varphi_{\eta})}^{-\hat{p}}\mathbb{E}\Big(\sup_{t\in[t_k,t_{k+1}]}\big|\int_{t_k}^{t} (t-u)^{-\eta}\big(\int_{0}^{u}(u-s)^{\eta-1} \beta\big(\mathbb{L}_{Z ^{N}(s)},\mathbb{L}_{Z ^{N}(s-\rho)}\big)        dB_{s}^{H,i}\big)du\big|^{\hat{p}}\Big)\\
  		\leq&{(\varphi_{\eta})}^{-\hat{p}}\mathbb{E}\Big(\sup_{t\in[t_k,t_{k+1}]}    \big[\big(\int_{t_k}^{t}(t-u)^{-\frac{\hat{p}\eta}{\hat{p}-1}}du\big)^{\hat{p}-1}   \big(\int_{t_k}^{t} \Big|\int_{0}^{u}(u-s)^{\eta-1} \beta\big(\mathbb{L}_{Z ^{N}(s)},\mathbb{L}_{Z ^{N}(s-\rho)}\big) dB_{s}^{H,i}   \Big| ^{\hat{p}}du  \big)
  		\big]                     
  		\Big)\\
  			\leq&\frac{{(\varphi_{\eta})}^{-\hat{p}}(\hat{p}-1)}{(\hat{p}-1-\hat{p}\eta)^{\hat{p}-1}}\Delta^{\hat{p}-1-\hat{p}\eta}\int_{t_k}^{t_{k+1}}\mathbb{E}\big|\int_{0}^{u}    (u-s)^{\eta-1} \beta\big(\mathbb{L}_{Z ^{N}(s)},\mathbb{L}_{Z ^{N}(s-\rho)}\big) dB_{s}^{H,i} \big|^{\hat{p}} du.
  	\end{split}
  \end{equation*}
By Theorem 1.1 in \cite{ccc27}, we see
  \begin{equation*}
	\begin{split}
	\mathbb{E}\big|\int_{0}^{u}    (u-s)^{\eta-1} \beta\big(\mathbb{L}_{Z ^{N}(s)},\mathbb{L}_{Z ^{N}(s-\rho)}\big) dB_{s}^{H,i} \big|^{\hat{p}} \leq C_{\hat{p},H}\big[      \int_{0}^{u}    (u-s)^\frac{\eta-1}{H}\big|\beta\big(\mathbb{L}_{Z ^{N}(s)},\mathbb{L}_{Z ^{N}(s-\rho)}\big) \big|^\frac{1}{H}ds
   \big]^{\hat{p}H}.
	\end{split}
\end{equation*}
By choosing $\tilde{q}=\hat{p}H$, $\alpha=\frac{H-1-\eta}{H}$, $\tilde{p}=\frac{\hat{p}H}{\hat{p}(\eta+H-1)+1}$ in Lemma 3.2 in \cite{ccc12}, we derive from Assumption \ref{cccass1} and Lemma \ref{cccnumbound} that
 \begin{equation*}
	\begin{split}
		&\mathbb{E}\big(\sup_{t\in[t_k,t_{k+1}]}|J_2(t)|^{\hat{p}}\big)\\
		\leq&\frac{{(\varphi_{\eta})}^{-\hat{p}}(\hat{p}-1)}{(\hat{p}-1-\hat{p}\eta)^{\hat{p}-1}}\Delta^{\hat{p}-1-\hat{p}\eta}\mathbb{E}\Big[\int_{t_k}^{t_{k+1}}  \Big| \beta\big(\mathbb{L}_{Z ^{N}(s)},\mathbb{L}_{Z ^{N}(s-\rho)}\big)\Big|^\frac{\hat{p}}{\hat{p}(\eta+H-1)+1} ds \Big]^{\hat{p}(\eta+H-1)+1} \\
		\leq&C_{\hat{p},H,\eta}\Delta^{\hat{p}-1-\hat{p}\eta}\Delta^{\hat{p}(\eta+H-1)}\int_{t_k}^{t_{k+1}} \mathbb{E} \big| \beta\big(\mathbb{L}_{Z ^{N}(s)},\mathbb{L}_{Z ^{N}(s-\rho)}\big)\big|^{\hat{p}}ds  \\
		\leq&C_{\hat{p},H,\eta,K_2}\Delta^{\hat{p}H-1}\int_{t_k}^{t_{k+1}} \mathbb{E} \big[ 1+\mathcal{W}_{2}^{\hat{p}}(\mathbb{L}_{Z ^{N}(s)})+\mathcal{W}_{2}^{\hat{p}}(\mathbb{L}_{Z ^{N}(s-\rho)})\big]ds  \\
		\leq&C_{\hat{p},H,\eta}\Delta^{\hat{p}H-1}\int_{t_k}^{t_{k+1}}  \Big( \mathbb{E}|\tilde{Z }^{i,N}(s)|^{\hat{p}}+\mathbb{E}|\tilde{Z }^{i,N}(s-\rho)|^{\hat{p}}\Big)ds  \\
		\leq&C_{\hat{p},H,\eta}\Delta^{\hat{p}H}.
	\end{split}
\end{equation*}
For $J_1(t)$, the H\"older inequality means that
 \begin{equation*}
	\begin{split}
		&\mathbb{E}\big(\sup_{t\in[t_k,t_{k+1}]}|J_1(t)|^{\hat{p}}\big)\\
		\leq&\Delta^{\hat{p}-1}\mathbb{E}\int_{t_k}^{t_{k+1}}  \big| \alpha\big(Z ^{i,N}(s),Z ^{i,N}(s-\rho),\mathbb{L}_{Z ^{N}(s)},\mathbb{L}_{Z ^{N}(s-\rho)}\big)\big|^{\hat{p}}ds  \\
		\leq&C_{\hat{p},K_1}\Delta^{\hat{p}-1}\mathbb{E}\int_{t_k}^{t_{k+1}}  \big[ 1+|Z ^{i,N}(s)|^{\hat{p}}+|Z ^{i,N}(s-\rho)|^{\hat{p}(l+1)}+\mathcal{W}_{2}^{\hat{p}}(\mathbb{L}_{Z ^{N}(s)})+\mathcal{W}_{2}^{\hat{p}}(\mathbb{L}_{Z ^{N}(s-\rho)})\big]ds  \\
		\leq&C_{\hat{p},\bar{p},T,H,K_1,K_2,\|\xi\|,\bar{p}_m}\Delta^{\hat{p}},
	\end{split}
\end{equation*}
  where Lemma \ref{cccnumbound} with $\hat{p}(l+1)\leq\bar{p}$ has been used. 
  Thus,
   \begin{equation*}
  	\begin{split}
  		\mathbb{E}\big(\sup_{t\in[t_k,t_{k+1}]}|\tilde{Z }^{i,N}(t)-Z ^{i,N}(t)|^{\hat{p}}\big)\leq&2^{\hat{p}-1}\mathbb{E}\big(\sup_{t\in[t_k,t_{k+1}]}|J_1(t)|^{\hat{p}}\big)+2^{\hat{p}-1}\mathbb{E}\big(\sup_{t\in[t_k,t_{k+1}]}|J_2(t)|^{\hat{p}}\big)\\
  			\leq&C_{\hat{p},\bar{p},T,H,K_1,K_2,\|\xi\|,\bar{p}_m}\Delta^{\hat{p}H}.
  	\end{split}
  \end{equation*}
 So the desired assertion holds. 
  \eproof
  
  	\begin{assp}\label{cccass2}
  	There exist constants $K_3>0$ and $\vartheta\in(0,1]$ such that, for any $\check{p}>0$,
  	\begin{equation*}
  		\mathbb{E}\left(\sup_{t,s\in[-\rho,0]}|\xi(t)-\xi(s)|^{\check{p}}\right)\leq K_3|t-s|^{\vartheta \check{p}}.
  	\end{equation*}
  \end{assp}
  
The following theorem reveals the convergence rate of EM scheme when $H\in (1/2,1)$.

  \begin{theorem}\label{ccctheo1}
  	For $H\in  (1/2,1)$, let Assumptions  \ref{cccass1} and \ref{cccass2} hold. For any $i\in\mathbb{S}_N$ and $p\in\left[2,\frac{\bar{p}}{2l}\right]$ with $\bar{p}>4l$, we see
  	\begin{equation*}
  		\begin{split}
  		 \mathbb{E}\big(\sup_{t\in[0,T]}|\tilde{Z }^{i,N}(t)-\Upsilon^{i,N}(t)|^p\big) \leq C_{p,\bar{p},T,H,K_1,K_2,K_3,\|\xi\|,\bar{p}_m,p_m}\Delta^{(\vartheta\wedge H)p}.\
  		\end{split}
  	\end{equation*} 
  \end{theorem}

  \noindent
  {\it Proof}.
   For any $t\in[0,T]$ and $i\in\mathbb{S}_N$,
applying  H\"older's inequality, Assumption \ref{cccass1}, (\ref{cccbbdt}) on (\ref{cccIPS}) and (\ref{ccccon1}) gives that
  \begin{equation*}
  	\begin{split}
  		&\mathbb{E}\big(\sup_{s\in[0,t]}|\tilde{Z }^{i,N}(s)-\Upsilon^{i,N}(s)|^p\big)\\
  		\leq&2^{p-1}\mathbb{E}\Big(\sup_{s\in[0,t]}\Big|\int_{0}^{s}\big[\alpha\big(Z ^{i,N}(u),Z ^{i,N}(u-\rho),\mathbb{L}_{Z ^{N}(u)},\mathbb{L}_{Z ^{N}(u-\rho)}\big)\\&-\alpha\big(\Upsilon^{i,N}(u),\Upsilon^{i,N}(u-\rho),\mathbb{L}_{\Upsilon^{N}(u)},\mathbb{L}_{\Upsilon^{N}(u-\rho)}\big)\big]du\Big|^p\Big)\\
  		&+2^{p-1}\mathbb{E}\Big(\sup_{s\in[0,t]}\Big|\int_{0}^{s}\big[\beta\big(\mathbb{L}_{Z ^{N}(u)},\mathbb{L}_{Z ^{N}(u-\rho)}\big)-\beta\big(\mathbb{L}_{\Upsilon^{N}(u)},\mathbb{L}_{\Upsilon^{N}(u-\rho)}\big)\big]dB_u^{H,i}\Big|^p\Big)\\
  		\leq&C_{p,T,H,K_1,K_2}\mathbb{E}\int_{0}^{t}\Big[|Z ^{i,N}(u)-\Upsilon^{i,N}(u)|^p\\&+\big(1+|Z ^{i,N}(u-\rho)|^l+|\Upsilon^{i,N}(u-\rho)|^l\big)^p|Z ^{i,N}(u-\rho)-\Upsilon^{i,N}(u-\rho)|^p\\&+\mathbb{W}_{p}^{p}(\mathbb{L}_{Z ^{N}(u)},\mathbb{L}_{\Upsilon^{N}(u)})+\mathbb{W}_{p}^{p}(\mathbb{L}_{Z ^{N}(u-\rho)},\mathbb{L}_{Z ^{N}(u-\rho)})\Big]du.
  	\end{split}
  \end{equation*}
  Notice that
  \begin{equation*}
  	\begin{split}
  		&\mathbb{E}\big[\big(1+|Z ^{i,N}(u-\rho)|^l+|\Upsilon^{i,N}(u-\rho)|^l\big)^p|Z ^{i,N}(u-\rho)-\Upsilon^{i,N}(u-\rho)|^p\big]\\
  		\leq&C_p\big[\mathbb{E}\big(1+|Z ^{i,N}(u-\rho)|^{2lp}+|\Upsilon^{i,N}(u-\rho)|^{2lp}\big)\big]^\frac{1}{2} \big[\mathbb{E}|Z ^{i,N}(u-\rho)-\Upsilon^{i,N}(u-\rho)|^{2p}\big]^\frac{1}{2}\\
  		\leq&C_{p,\bar{p},T,H,K_1,K_2,\|\xi\|,\bar{p}_m}\Big(\big[\mathbb{E}|Z ^{i,N}(u-\rho)-\tilde{Z }^{i,N}(u-\rho)|^{2p}\big]^\frac{1}{2}+\big[\mathbb{E}|\tilde{Z }^{i,N}(u-\rho)-\Upsilon^{i,N}(u-\rho)|^{2p}\big]^\frac{1}{2}\Big),
  	\end{split}
  \end{equation*}
where Lemma \ref{cccnumbound} with $2lp\leq\bar{p}$ is used. Then by
$
  	\mathbb{E}\mathbb{W}_{p}^{p}(\mathbb{L}_{Z ^{N}(\cdot)},\mathbb{L}_{\Upsilon^{N}(\cdot)})
  		\leq \frac{1}{N}\sum_{j=1}^{N}\mathbb{E}|Z ^{j,N}(\cdot)-\Upsilon^{j,N}(\cdot)|^{p},
  $
we derive from Lemma \ref{ccccondis} that
    \begin{equation*}
  	\begin{split}
  		&\mathbb{E}\big(\sup_{s\in[0,t]}|\tilde{Z }^{i,N}(s)-\Upsilon^{i,N}(s)|^p\big)\\
  		\leq&C_{p,\bar{p},T,H,K_1,K_2,\|\xi\|,\bar{p}_m}\int_{0}^{t}\Big[ \mathbb{E}\big(\sup_{u\in[0,s]}|Z ^{i,N}(u)-\tilde{Z }^{i,N}(u)|^p\big)+\mathbb{E}\big(\sup_{u\in[0,s]}|\tilde{Z }^{i,N}(u)-\Upsilon^{i,N}(u)|^p\big)\\&+\mathbb{E}\big(\sup_{u\in[-\rho,0]}|\xi(u)-\xi(\lfloor\frac{u}{\Delta}\rfloor\Delta)|^p\big)+\big(\mathbb{E}|\tilde{Z }^{i,N}(s-\rho)-\Upsilon^{i,N}(s-\rho)|^{2p}\big)^\frac{1}{2}+\Delta^{pH}\Big]ds.
  	\end{split}
  \end{equation*}
The Gronwall inequality with Assumption \ref{cccass2} and Lemma \ref{ccccondis} leads to
  \begin{equation}\label{cccconver1}
	\begin{split}
		&\mathbb{E}\big(\sup_{s\in[0,t]}|\tilde{Z }^{i,N}(s)-\Upsilon^{i,N}(s)|^p\big)\\
		\leq&C_{p,\bar{p},T,H,K_1,K_2,K_3,\|\xi\|,\bar{p}_m}\Big(\Delta^{(\vartheta\wedge H)p}+\big[\mathbb{E}\big(\sup_{s\in[0,t]}|\tilde{Z }^{i,N}(s-\rho)-\Upsilon^{i,N}(s-\rho)|^{2p}\big)\big]^\frac{1}{2}\Big).
	\end{split}
\end{equation}
Define a sequence as 
\begin{equation*}
p_{m}=(2-m+\lfloor\frac{T}{\rho}\rfloor)p2^{1-m+\lfloor\frac{T}{\rho}\rfloor},~~~~m=1, 2, \cdots, \lfloor\frac{T}{\rho}\rfloor+1.
\end{equation*}
One can see $p_{\lfloor\frac{T}{\rho}\rfloor+1}=p$ and $2 p_{m+1}<p_{m}$  for $m=1, 2, \cdots, \lfloor\frac{T}{\rho}\rfloor$.
For $s\in[0,\rho]$, we get from (\ref{cccconver1}) that
\begin{equation}\label{cccconver2}
	\begin{split}
		\mathbb{E}\big(\sup_{s\in[0,\rho]}|\tilde{Z }^{i,N}(s)-\Upsilon^{i,N}(s)|^{p_1}\big)
		\leq&C_{p_1,\bar{p},T,H,K_1,K_2,K_3,\|\xi\|,\bar{p}_m}\Delta^{(\vartheta\wedge H)p_1}.
	\end{split}
\end{equation}
For $s\in[0,2\rho]$, the H\"older inequality with (\ref{cccconver1}) and (\ref{cccconver2}) yields that
\begin{equation*}
	\begin{split}
		&\mathbb{E}\big(\sup_{s\in[0,2\rho]}|\tilde{Z }^{i,N}(s)-\Upsilon^{i,N}(s)|^{p_2}\big)\\
		\leq&C_{p_1,p_2,\bar{p},T,H,K_1,K_2,K_3,\|\xi\|,\bar{p}_m}\Big(\Delta^{(\vartheta\wedge H)p_2}+\big[\mathbb{E}\big(\sup_{s\in[0,2\rho]}|\tilde{Z }^{i,N}(s-\rho)-\Upsilon^{i,N}(s-\rho)|^{2p_2}\big)\big]^\frac{1}{2}\Big)\\
		\leq&C_{p_1,p_2,\bar{p},T,H,K_1,K_2,K_3,\|\xi\|,\bar{p}_m}\Big(\Delta^{(\vartheta\wedge H)p_2}+\big[\mathbb{E}\big(\sup_{s\in[0,2\rho]}|\tilde{Z }^{i,N}(s-\rho)-\Upsilon^{i,N}(s-\rho)|^{p_1}\big)\big]^\frac{p_2}{p_1}\Big)\\
			\leq&C_{p_1,p_2,\bar{p},T,H,K_1,K_2,K_3,\|\xi\|,\bar{p}_m}\Delta^{(\vartheta\wedge H)p_2}.
	\end{split}
\end{equation*}
The induction about the time segment generated by time-delay $\rho$ gives that
 \begin{equation*}
	\begin{split}
		\mathbb{E}\big(\sup_{t\in[0,T]}|\tilde{Z }^{i,N}(t)-\Upsilon^{i,N}(t)|^p\big)
		\leq&C_{p,\bar{p},T,H,K_1,K_2,K_3,\|\xi\|,\bar{p}_m,p_m}\Delta^{(\vartheta\wedge H)p},
	\end{split}
\end{equation*}
where $m=1,2 \cdots,\lfloor\frac{T}{\rho}\rfloor+1$.
\eproof

As for the convergence rate of EM scheme in the case $H\in (0,1/2)$, the Corollary 5.5.3 in \cite{ccc2} plays a key role, so we quote it as the following lemma.
\begin{lemma}\label{ccchxiao}
	For each $\check{p} >1$, there is a $C_{\check{p}}<\infty$ such that
	\begin{equation*}
		\begin{split}
		 \mathbb{E}|B_t^H-B_s^H|^{\check{p}}\leq C_{\check{p}}|t-s|^{\check{p}H}.
		\end{split}
	\end{equation*} 
\end{lemma}

\begin{lemma}\label{cccd2d1}
	For $H\in  (0,1/2)$, let Assumption \ref{cccass1} hold. For any $i\in\mathbb{S}_N$, $t\in[0,T]$ and $\hat{p}\in\left[2,\frac{\bar{p}}{l+1}\right]$ with $\bar{p}>2l+2$, we have
	\begin{equation*}
		\begin{split}
			 \mathbb{E}|\tilde{Z }^{i,N}(t)-Z ^{i,N}(t)|^{\hat{p}} \leq C_{\hat{p},\bar{p},T,H,K_1,\beta,\|\xi\|,\bar{p}_m}\Delta^{\hat{p}H}.
		\end{split}
	\end{equation*} 
\end{lemma}
{\it Proof}.
For any $t\in[t_k,t_{k+1})$ and $i\in\mathbb{S}_N$, we see
\begin{equation*}
	\begin{split}
		&\mathbb{E}|\tilde{Z }^{i,N}(t)-Z ^{i,N}(t)|^{\hat{p}}\\\leq&2^{\hat{p}-1}\mathbb{E}\left|\int_{t_k}^{t}\alpha\big(Z ^{i,N}(s),Z ^{i,N}(s-\rho),\mathbb{L}_{Z ^{N}(s)},\mathbb{L}_{Z ^{N}(s-\rho)}\big)ds\right|^{\hat{p}}+2^{\hat{p}-1}\mathbb{E}\left|\int_{t_k}^{t}\beta dB_{s}^{H,i}\right|^{\hat{p}}\\
		\leq & C_{\hat{p},\bar{p},T,H,K_1,\beta,\|\xi\|,\bar{p}_m}\Delta^{\hat{p}H},
	\end{split}
\end{equation*}
where we have used Lemma \ref{ccchxiao} and the estimation of $J_1$ in the proof of Lemma \ref{ccccondis}.
\eproof

 \begin{theorem}\label{ccctheo2}
	For $H\in  (0,1/2)$, let Assumptions  \ref{cccass1} and \ref{cccass2} hold. For any $i\in\mathbb{S}_N$ and $p\in\left[2,\frac{\bar{p}}{2l}\right]$ with $\bar{p}>4l$, we have
	\begin{equation*}
		\begin{split}
			\mathbb{E}|\tilde{Z }^{i,N}(T)-\Upsilon^{i,N}(T)|^p \leq C_{p,\bar{p},T,H,K_1,\beta,K_3,\|\xi\|,\bar{p}_m,p_m}\Delta^{(\vartheta\wedge H)p}.\
		\end{split}
	\end{equation*} 
\end{theorem}
Thanks to Lemma \ref{cccd2d1}, the proof of Theorem \ref{ccctheo2} is similar to Theorem \ref{ccctheo1}, since the stochastic integral  vanishes in $\tilde{Z }^{i,N}(t)-\Upsilon^{i,N}(t)$.
Then the following main conclusion is reached by exploiting the triangle inequality for Theorems \ref{cccfifi1}, \ref{ccctheo1} and \ref{ccctheo2}.
\begin{theorem}
	For $H\in  (1/2,1)$, let all conditions in Theorems \ref{cccfifi1} and \ref{ccctheo1} be satisfied. Then we have
\begin{equation*}
	\begin{split}
		\mathbb{E}\Big(\sup_{t\in[0,T]}|\Upsilon^{i}(t)-\tilde{Z }^{i,N}(t)|^{p}\Big)\leq C
		\left\{\begin{array}{ll}
			(N^{-1 / 2})^{\lambda_{p,T,\rho}}+\Delta^{(\vartheta\wedge H)p}, & \text { if } p>d/2, \\
			{[N^{-1 / 2} \log(1+N)]}^{\lambda_{p,T,\rho}}+\Delta^{(\vartheta\wedge H)p}, & \text { if } p=d/2,\\
			{(N^{-p / d})}^{\lambda_{p,T,\rho}}+\Delta^{(\vartheta\wedge H)p}, & \text { if }2\leq p<d/2,
		\end{array}\right.
	\end{split}
\end{equation*}
where $\lambda_{p,T,\rho}=(\frac{p-\varepsilon}{p})^{\lfloor\frac{T}{\rho}\rfloor}$ and $C$ is a positive real constant dependent of $p,\bar{p}, $, $T$, $H$, $K_1$, $K_2$, $K_3$,$\|\xi\|$, $ \bar{p}_m$,$ p_m$ but independent of $N$, $\Delta$.

For $H\in  (0,1/2)$, let all conditions in Theorems \ref{cccfifi1} and \ref{ccctheo2} be satisfied and the diffusion coefficient be a constant $\beta$. Then we have
	\begin{equation*}
		\begin{split}
			\mathbb{E}|\Upsilon^{i}(T)-\tilde{Z }^{i,N}(T)|^{p}\leq C
			\left\{\begin{array}{ll}
				(N^{-1 / 2})^{\lambda_{p,T,\rho}}+\Delta^{(\vartheta\wedge H)p}, & \text { if } p>d/2, \\
				{[N^{-1 / 2} \log(1+N)]}^{\lambda_{p,T,\rho}}+\Delta^{(\vartheta\wedge H)p}, & \text { if } p=d/2,\\
				{(N^{-p / d})}^{\lambda_{p,T,\rho}}+\Delta^{(\vartheta\wedge H)p}, & \text { if }2\leq p<d/2,
			\end{array}\right.
		\end{split}
	\end{equation*}
where $\lambda_{p,T,\rho}=(\frac{p-\varepsilon}{p})^{\lfloor\frac{T}{\rho}\rfloor}$ and $C$ is a positive real constant dependent of $p,\bar{p}, $, $T$, $H$, $K_1$,  $\beta$, $K_3$, $\|\xi\|$, $ \bar{p}_m$,$ p_m$ but independent of $N$, $\Delta$.

\end{theorem}

\section{Numerical experiments}\label{cccsec05}

\begin{figure}
	\centering
	\subfigure[H=0.6]{
	\includegraphics[width =0.38\linewidth]{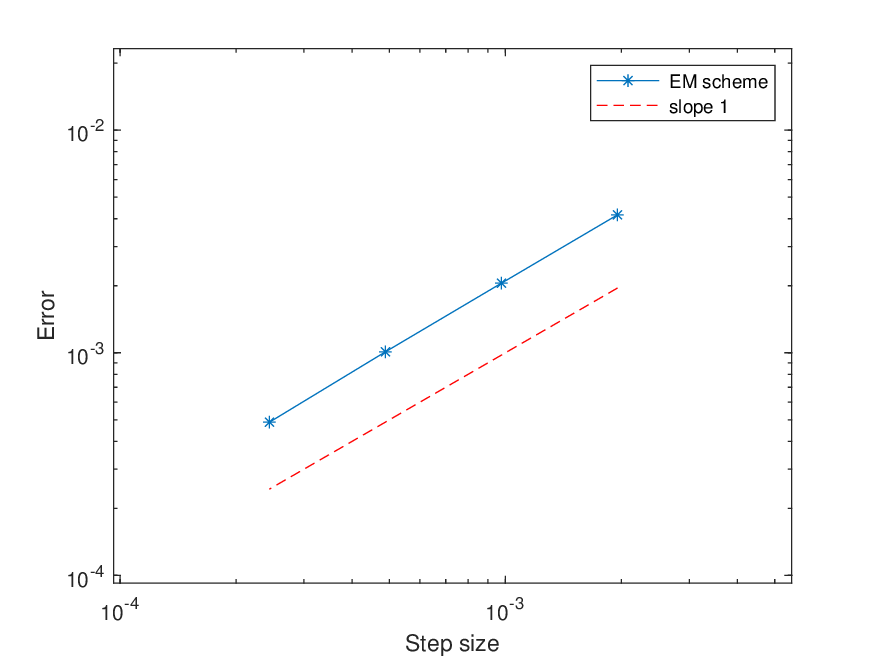}
}
\subfigure[H=0.7]{
	\includegraphics[width =0.38\linewidth]{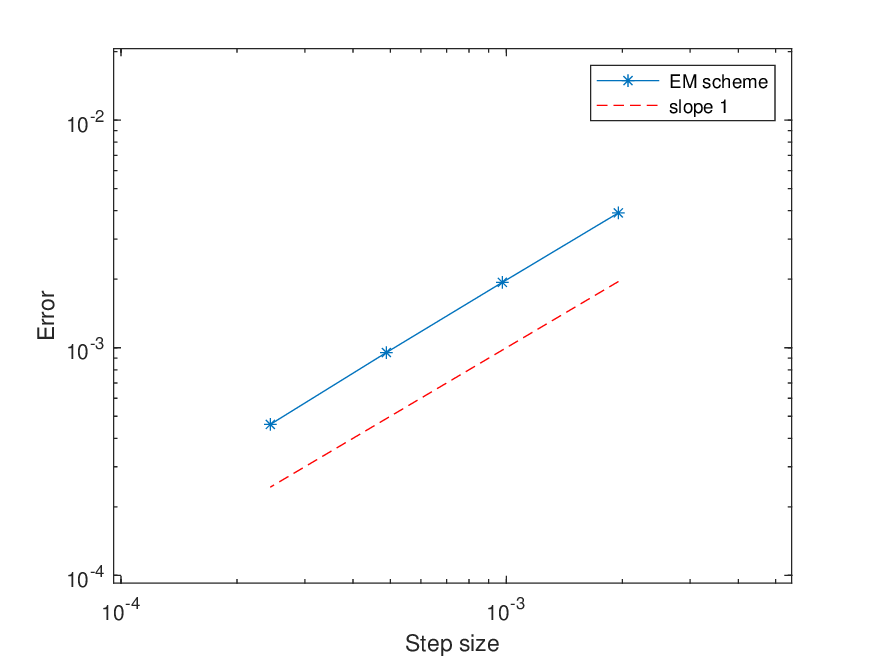}
}\\
\subfigure[H=0.8]{
	\includegraphics[width =0.38\linewidth]{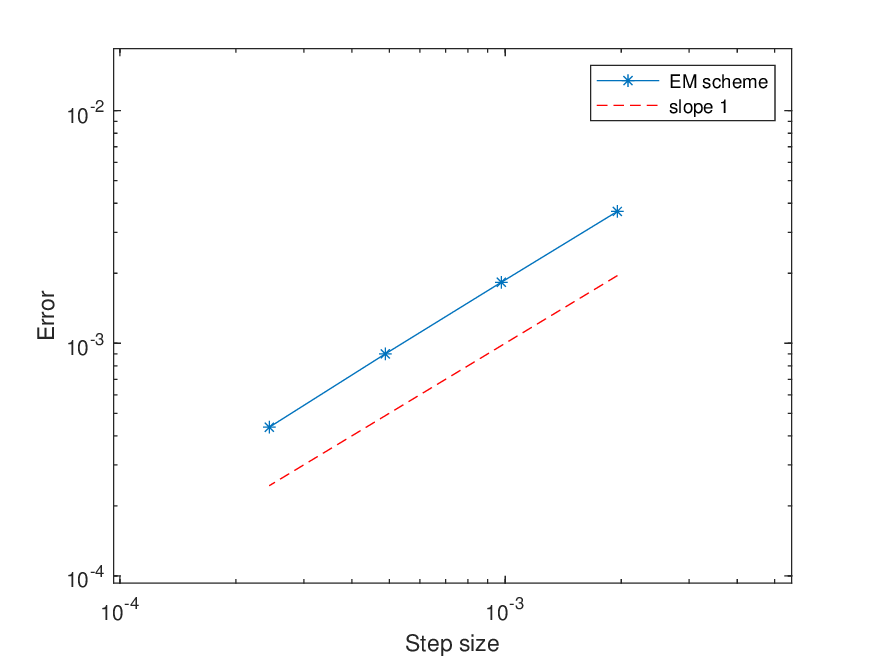}
}
\subfigure[H=0.9]{
	\includegraphics[width =0.38\linewidth]{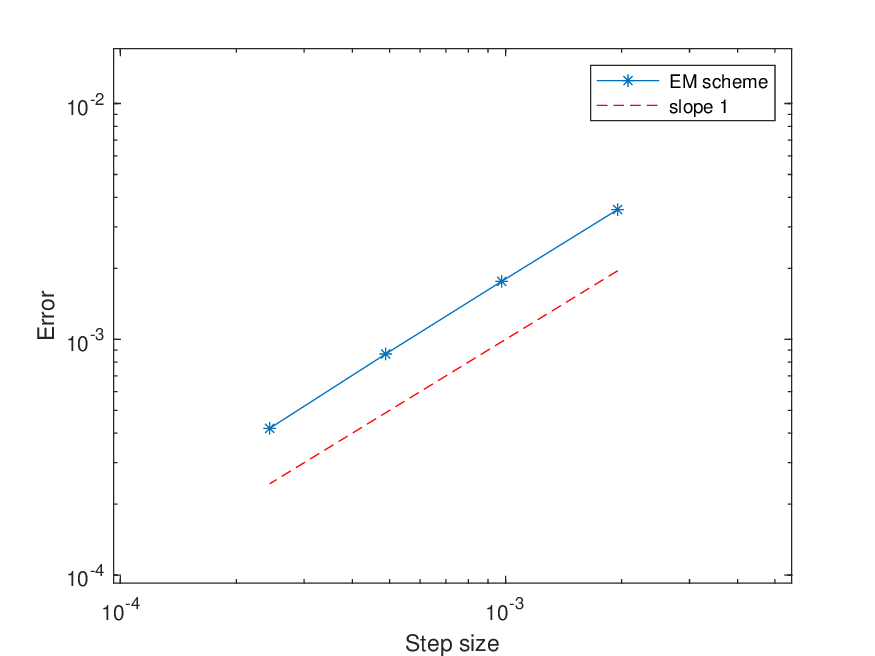}
}
	\caption{Convergence rates of classical EM scheme for (\ref{cccexa1})\label{cccffig}}
\end{figure}

Based on the statement of SODM in Section \ref{cccsec01}, we now consider the following scalar equation
\begin{equation} \label{cccexa1}
	\begin{split}
d\Upsilon(t) =&\big[a_1\int_{\mathbb{R}}\Phi(|\Upsilon(t)-x|)(\Upsilon(t)-x)\mathbb{L}_\Upsilon (dx)+a_2\Upsilon(t)+a_3\psi(\Upsilon(t-\rho))+a_4\mathbb{E}\Upsilon(t-\rho)\big]dt \\
&+a_5 dB_t^H, \quad t\in[0,T],	
\end{split}
\end{equation}
with the initial value $\xi(\theta)=|\theta|, \theta \in [-\rho,0]$, where the interaction kernel $\Phi(\cdot)$ represents the affect of an individual
to each other, and $\psi(\Upsilon(t-\rho))$ is the function of $\Upsilon(t-\rho)$.
Obviously, this stochastic
opinion dynamics model can reflect the effects of both intrinsic memory (delay state variable and its law) and extrinsic memory (fracional noise).
For more details about such SODM, please refer to \cite{ccc32,ccc33,ccc34}.
Let $a_1=a_2=a_4=a_5=1$, $a_3=-1$, $\rho=\frac{1}{8}$ and $\psi(\Upsilon(t-\rho))=\Upsilon(t-\rho)^3$,
\begin{equation*}
	\begin{split}
\Phi(|\Upsilon(t)-x|)&=\begin{cases}\sin(|\Upsilon(t)-x|-0.5),& \text{if}~~ |\Upsilon(t)-x|<0.5,\\
	\cos(|\Upsilon(t)-x|+0.5),&\text{if}~~ |\Upsilon(t)-x|\geq 0.5.\end{cases}
\end{split}
\end{equation*}
The assumptions are satisfied. We use the numerical solution with $\Delta=2^{-16}$ as the analytical solution since it cannot be represented explicitly.
Then the numeircal simulation 
for the error which is defined by
	\begin{equation}\label{cccerer1}
	err=\left[\frac{1}{N}\sum_{i=1}^{N}|\Upsilon_{(\varsigma)}^{i,N}(T)-Z _{(\varsigma_r)}^{i,N}(T)|^2\right]^{\frac{1}{2}},
\end{equation}	
with $T=1,N=200,\rho=\frac{1}{8}$
is implemented, where $\varsigma$ and $\varsigma_r (r\in\{1,2,3,4\})$ mean the level of the time discretization. Here, $\Upsilon_{(\varsigma)}^{i,N}(T)$ is regarded as the solution to interacting particle system (\ref{cccIPS}) with $\Delta=2^{-16}$ at $T=1$, and $Z _{(\varsigma_r)}^{i,N}(T)$ is the numerical solution to classical EM scheme (\ref{cccemli}) at $T=1$
with $\varsigma_r \in \{\varsigma_1,\varsigma_2,\varsigma_3,\varsigma_4\}$ matching $\Delta\in\{2^{-9}, 2^{-10},2^{-11}, 2^{-12}\}$.
  The errors between interacting particle system  and numerical solution are depicted in Figure \ref{cccffig},
 from which we can observe that the slope of these lines can support the theoretical results.
 
\section*{Funding}

This work is supported by the National Natural Science Foundation of China (12271368, 62373383 and 62076106) and Fund for Academic Innovation Teams of South-Central Minzu University
(XTZ24004).

\end{document}